\documentclass[12pt,reqno]{amsart}
\usepackage{mathrsfs}
\usepackage{geometry}
\usepackage{lineno}
\usepackage[all]{xy}
\usepackage{tikz}
\usetikzlibrary{positioning}
\tikzset{>=stealth}

\usepackage{enumerate} 
\usepackage[utf8]{inputenc} 
\usepackage[T1]{fontenc} 
\usepackage{amssymb,amsthm}
\usepackage{mathrsfs}
\usepackage{textcomp}
\usepackage{bm}
\usepackage{graphicx}
\usepackage{xcolor}

\newtheorem{thm}{Theorem}[section]
\newtheorem{prop}[thm]{Proposition}
\newtheorem{lem}[thm]{Lemma}
\newtheorem{cor}[thm]{Corollary}

\theoremstyle{definition}

\theoremstyle{remark}

\numberwithin{equation}{section}

\usepackage[font=footnotesize,tableposition=top,figureposition=bottom,skip=0pt]{caption}
\usepackage{hyperref}
\hypersetup{%
	linkcolor=blue,
	anchorcolor=black,
	citecolor=red,
	filecolor=magenta,
	menucolor=red,
	urlcolor=magenta,
	colorlinks=true,
	bookmarksnumbered=true,
	bookmarksopen=true,
	bookmarksopenlevel=0,
	breaklinks=true
}
\usepackage{color}

\newcommand{\N}{\mathbb{N}}  
\newcommand{\Z}{\mathbb{Z}} 

\begin{document}


\title
{ Bohr chaoticity of topological dynamical  systems}

 \author{Aihua Fan}
\address{(A. H. FAN) 
	LAMFA, UMR 7352 CNRS, University of Picardie, 33 rue Saint Leu,80039 Amiens, France}
\email{ai-hua.fan@u-picardie.fr}


\author{Shilei Fan}
\address{(S. L. Fan) School of Mathematics and Statistics, Hubei Key Laboratory of Mathematical Sciences, Central  China Normal University,  Wuhan, 430079, China}
\email{slfan@mail.ccnu.edu.cn}

\author{Valery V. Ryzhikov}
\address{(V. V. Ryzhikov)
Moscow State University Lomonosov, Faculty of Mechanics and Mathematics, Moscow, Russia, 119991}
\email{vryzh@mail.ru}

\author{Weixiao Shen}
\address{(W. X. Shen) Shanghai Center for Mathematical Sciences, Fudan University, 220 Handan Road, Shanghai 200433, China }
\email{ wxshen@fudan.edu.cn}


\thanks{A. H. FAN  was supported by NSF of China (Grant No. 11971192); S. L. FAN was supported by NSF of China (Grant No. 11971190) and Fok Ying-Tong Education Foundation, China (Grant No.171001);
W. X. SHEN was supported by NSF of China (Grant No 11731003). }
\maketitle




\begin{abstract}
	We introduce the notion of Bohr chaoticity, which is a topological invariant for topological dynamical systems, and which is opposite to the property required by Sarnak's conjecture. We prove the Bohr chaoticity for all systems which have a horseshoe and
	for all toral affine dynamical systems  of positive entropy, some of which don't have a horseshoe.
	But uniquely ergodic dynamical systems are not Bohr chaotic.
\end{abstract}

\section{Introduction}
A sequence of complex numbers  $(w_n)_{n\ge 0}\in \ell^\infty(\mathbb{N})$  is called a (non-trivial) {\em weight}
or {\em weight sequence}  if it satisfies
\begin{equation}\label{w}
\limsup_{N\to \infty} \frac{1}{N} \sum_{n=0}^{N-1}|w_n|>0.
\end{equation}
Let us consider a weight $(w_n)$ and a topological dynamical system  $(X, T)$. For a given continuous function  $f\in C(X)$ and a given point $x\in X$, we say that $(w_n)$ is {\em orthogonal} to the observation $(f(T^n x))_{n\ge 0}$ if
\begin{equation}\label{Orth}
\lim_{N\to \infty} \frac{1}{N} \sum_{n=0}^{N-1}w_n f(T^n x)=0.
\end{equation}
If it is the case for all $f\in C(X)$ and all $x\in X$, we say that
$(w_n)$ is orthogonal to the dynamical system $(X, T)$.
	Sarnak's conjecture
states that the M\"{o}bius sequence $(\mu(n))$ is orthogonal to all
topological dynamical systems of zero entropy \cite{Sarnak}. This conjecture is still open. But it is known that there do exist non-trivial weights other than the M\"{o}bius sequence
which are orthogonal to all topological dynamical systems of zero entropy, for example almost all symmetric Bernoulli sequence taking
values $-1$ and $1$ (see \cite{AKLD2017}). It is also the case for $(e^{2\pi i \beta^n})$ for almost all $\beta >1$ \cite{Shi}.

We consider a  property of dynamical systems that is completely opposite to that property required by Sarnak's conjecture.
A topological dynamical system  $(X, T)$ is said to be {\em Bohr chaotic} if it is not orthogonal to any (non-trivial) weights. That is to say,  for any non-trivial weight $(w_n)_{n\ge 0} \in \ell^\infty(\N)$,
there exist a continuous function $f \in C(X)$ and a point $x\in X$ such that
\begin{equation}\label{BC}
\limsup_{N\to \infty} \frac{1}{N} \left|\sum_{n=0}^{N-1}w_n f(T^n x)\right|>0.
\end{equation}
We say that $(f, x)$ verifying (\ref{BC}) is a  pair correlated to $(w_n)$ for the system
$(X,T)$.
\medskip

It is easy to see that the Bohr chaoticity is a topological invariant,
namely topologically conjugate systems share the Bohr chaoticity at the same time. It is even true that the extensions of a Bohr chaotic system are all Bohr chaotic (Proposition \ref{factor}).  In this paper we will provide
some Bohr chaotic systems and some non Bohr chaotic systems among systems of positive entropy.
\medskip

Bohr chaoticity is a kind of complexity, just like entropy positivity, Li-York chaoticity etc.
A Bohr chaotic system must have positive entropy, according to the preceding discussion. But having positive entropy is not sufficient for having Bohr chaoticity. In other words,  there are dynamical systems of positive entropy which are not Bohr-chaotic and it is the case for uniquely ergodic systems
(see Theorem \ref{thm:UE} below).
\medskip

Recall that a topological dynamical system $(X, T )$ admits a {\em two-sided horseshoe} (resp. {\em one-sided horseshoe}) if there exists a subsystem $(\Lambda, T^N)$ of some  power system $(X, T^N)$ ($N\ge 1$ being some integer), that is topologically conjugate  to the two-sided  full shift  $(\{0, 1\}^{\Z}, \sigma)$ (resp. one-sided full shift $(\{0, 1\}^{\N}, \sigma)$). The subsystem  $(\Lambda, T^N )$ is sometimes called a {\em $N$-order horseshoe} of the system $(X,T)$.
If, furthermore, the horseshoe $(\Lambda, T^N)$ satisfies the extra condition
\[\Lambda \cap \bigcup_{k=1}^{N-1} T^{k}(\Lambda)=\emptyset,\]
we say that the $N$-order horseshoe $(\Lambda, T^N )$ has {\em disjoint steps}, meaning that the steps $T^k(\Lambda)$ ($1\le k\le N-1$) are disjoint from the initial set $\Lambda$.

We shall prove that systems having a horseshoe must have a horseshoe with disjoint steps (Theorem \ref{thm:disjiontHS}).
The proof of this fact occupies an important part of the present paper. But it has a nice consequence on the Bohr chaoticity of the system.

\begin{thm}\label{thm:horse}
	Any system having a horseshoe is
	Bohr chaotic.
\end{thm}

As corollaries of Theorem \ref{thm:horse},
	all subshifts of finite type of positive entropy are Bohr-chaotic and
	all $\beta$-shifts are Bohr chaotic.
Actually all piecewise monotonic interval maps  of positive entropy are Bohr chaotic, because
Young \cite{Young1981} proved the existence of subsystems which are subshifts of finite type of positive entropy.
\medskip

By a result of Smale \cite{Smale1965}, on a manifold
$X$  of dimension $\dim X \ge 2$
there are many diffeomorphisms which
are Bohr Chaotic.
Let us state Smale's precise  result:
there exists a non-empty open set $U$, in the $C^1$-topology,   of the space of diffeomorphisms of $X$ such that for  each $T\in U$ there  exists  a $T$-invariant Cantor  set $K\subset X$ such that the subsystem $T: K\to K$ is topologically conjugate to a full shift. 	Consequently, by Proposition \ref{factor} and Proposition \ref{FS}, every $T\in U$ is Bohr chaotic. Here we don't need Theorem \ref{thm:horse}, because it is easy to see that the full shift is Bohr chaotic.
\medskip

	For an individual smooth system, let us state that every  $C^{1+\alpha}$ ($\alpha >0$)
diffeomorphism $T$ of a compact smooth manifold
admitting an ergodic non-atomic Borel probability invariant measure with non-zero Lyapunov exponents is Bohr chaotic.
 Because Katok \cite{Katok1980}  proved the existence of
a closed invariant hyperbolic set $F$ such that the subsystem $T: F\to F$ is topologically conjugate to a subshift of finite type (called
topological Markov chain in \cite{Katok1980}) of positive entropy.
In particular, Anosov systems are Bohr chaotic.
\medskip

But every Bohr chaotic system doesn't necessarily have a horseshoe. Such systems exist among toral automorphisms, as we shall see.
Let us recall that an affine map on the torus $\mathbb{T}^d$ ($d\ge 1$)
is of  the form $T x = Bx+b  \mod \mathbb{Z}^d$ where $b\in \mathbb{R}^d$ and  $B$ is a $d\times d$ matrix
of integral entries. We assume that $\det B \not=0$. By Sinai's theorem, the topological entropy $T$ is equal to $\sum_{i} \log |\lambda_i|$ where the sum is taken over all eigenvalues of $B$ such that $|\lambda_i|>1$ (cf. Theorem 8.15, \cite{Walters1982}).

\begin{thm}\label{torus}
	All toral affine systems of positive entropy are Bohr-chaotic.
	Actually the set of $x\in \mathbb{T}^d$  such that \eqref{BC} holds
	has Hausdorff dimension $d$.
\end{thm}

According to Lind and Schmidt (cf. Example 3.4 in \cite{LS1999}),
any ergodic irreducible partially hyperbolic toral automorphism
$T_A$ defined by a matrix $A\in {\rm GL}(d, \mathbb{Z})$ has no non-trivial homoclinic point.   Here is an example:
\begin{equation}\label{eq:AA}
   A= \begin{pmatrix}
      0 & 1 & 0 & 0\\
      0 & 0 & 1 & 0\\
      0 & 0 & 0 & 1\\
      -1 & 3 & -3 & 3
   \end{pmatrix}
\end{equation}
which is the companion matrix of the polynomial $z^4-3z^3 +3z^2 -3z +1$.
Since the fullshift
$(\{0,1\}^\mathbb{Z}, \sigma)$ has a dense set of homoclinic points, automorphisms like $T_A$ can not have horseshoe.
But it is Bohr chaotic by Theorem \ref{torus}.

\medskip

In a forthcoming paper \cite{FSV2019}, the authors have studied the Bohr chaoticity
for principal algebraic $\mathbb{Z}^d$-actions on $\mathbb{T}^{\mathbb{Z}^d}$. When $d>1$, a key point of proof is the existence of summable homoclinic points.
\medskip

Now we point out that the unique ergodicity is an obstruction for the Bohr chaoticity, as we state in the following theorem.

\begin{thm}\label{thm:UE}
	The uniquely ergodic dynamical systems are not Bohr chaotic.
\end{thm}

Theorem \ref{thm:UE} is a consequence of a result due to Robinson \cite{Rob1994} about the uniform convergence of Wiener-Wintner ergodic averages.  Recall that there exist uniquely ergodic dynamical systems having  positive entropy according to
Krieger \cite{Krieger1972}. See \cite{BCL2007} for smooth uniquely ergodic dynamical systems having positive entropy.
\medskip

Theorem \ref{torus} will be proved by using Riesz product measures borrowed from
harmonic analysis (cf. \cite{HZ1966,Zygmund2002}), in order to find points $x$ satisfying (\ref{BC}). As we will see in the proof, for any non-trivial weight,  the weighted ergodic
limit in \eqref{BC} exists and is equal to a  constant different from zero for
almost every point with respect to the Riesz product measure.
We should point out that these Riesz product measures are not $T$-invariant, because it is proved that the weighted ergodic limit in (\ref{BC}) is zero for any invariant measure and any oscillating weight of Davenport type (cf. \cite{F2017b}).  The positivity of entropy implies some lacunarity of the affine map (cf. Lemma \ref{lem:1} and Lemma \ref{lem2}) which effectively allows us to construct such measures, and the convergence of the weighted ergodic averages is proved by using the classical Menshov-Redemacher theorem. This method is also used in \cite{FSV2019}.
\medskip

Theorem \ref{thm:horse} will be proved by another method.  From the given horseshoe $(\Lambda, T^N)$,  we shall construct a $T^\tau$-invariant set $K$ for some possibly very large integer $\tau\ge 1$ such that $K$ is disjoint from the next steps
$\bigcup_{1\le j\le \tau-1} T^j(K)$ of $K$ under the action of $T$ before going back to $K$.  This system $(K, T^\tau)$ is Bohr chaotic and the Bohr chaoticity of
$(X, T)$ follows (cf. Theorem \ref{ThmReturn}). A little more can be proved. Actually the system $(K, T^\tau)$ that we construct  is a horseshoe with disjoint steps (cf. Theorem \ref{thm:disjiontHS}).
This is the main part of the proof: construct a horseshoe with disjoint steps from a horseshoe $(\Lambda, T^N)$. To do so, we consider the extended horseshoe $T: \Lambda^* \to \Lambda^*$
 where
$$
    \Lambda^* = \bigcup_{n=0}^{N-1} T^n(\Lambda).
$$
This is a subsystem of the original system $(X, T)$. We shall use the fact that the extended horseshoe doesn't have root, that is to say,  the equation $S^n=T$ has no continuous solution $S$ for any $n\ge 2$ (cf. Proposition \ref{prop:no-root}).
This fact forces that $\Lambda \not\subset \bigcup_{j\in J}T^j(\Lambda)$ for some set $J$ (cf. Proposition \ref{PO2}).
We say that $\Lambda$ is displaced by $T^j$ with $j\in J$.  Roughly speaking,
we can construct a better horseshoe in  $\Lambda \setminus \bigcup_{j\in J}T^j(\Lambda)$.  Here by ``better'' we means that the new horseshoe has ``more'' disjoint steps.
The possibility to construct this new horseshoe is based on the fact that in any cylinder of the full shift space we can construct a horseshoe with disjoint steps (cf. Proposition \ref{lem-SBSYS}).
If the order $N$ of the horseshoe is equal to $2$, this new horseshoe is already what we need. But for general order $N$, we have to repeat this procedure many times (cf. Proposition \ref{prop:disp1},
Proposition \ref{prop:disp2}, Lemma \ref{lem-final}).
\medskip

At the end of this introduction we mention two questions to which we don't have answer.
\medskip

Recall that {\em semi-horseshoes} are similarly defined as horseshoes are defined, but  the conjugation is weaken to the semi-conjugation.
Semi-horseshoes
are
well studied in \cite{HLXY2019,HY2011,KOR2016,LY2011}
\medskip.

{\bf Question 1.} {\it Are systems having semi-horseshoes Bohr-chaotic ?}
\medskip

Herman \cite{Herman1981} constructed a real-analytic diffeomorphism on a compact, connected manifold of dimension $4$ that is minimal and has positive entropy. Because of the minimality, Herman's diffeomorphism doesn't has horseshoe,  so Theorem \ref{thm:horse}  doesn't apply to it. 
Herman's diffeomorphism is not uniquely ergodic, so Theorem \ref{thm:UE}  doesn't apply either.
\medskip

{\bf Question 2.}  {\it Is Herman's diffeomorphism not Bohr-chaotic ?}

\medskip

\medskip

Finally let us explain the organization of the paper.
We shall first prove Theorem \ref{thm:UE} in Section \ref{Sect:UE}, the proof of  which is easy.
In Section \ref{Sect:subsystems} we shall study the Bohr chaoticity of a system
by that of its subsystems, Proposition \ref{factor} there shows that the Bohr chaoticity  is a topological invariant  and Theorem \ref{ThmReturn} there will serve as a tool
 to prove Theorem \ref{thm:horse}.
  The Bohr-chaoticity of affine maps on torus
  (Theorem \ref{torus}) will be proved in Section \ref{Sect:Affine}.  Theorem \ref{thm:horse} will be proved in Section \ref{sect:1horseshoe} (the case of one-sided horseshoe) and Section
  \ref{sect:2horseshoe} (the case of two-sided horseshoe).  We put two technical results in Appendixes, which could be known. Appendix A proves that any cylinder contains
  a horseshoe. Appendix B provides two proofs of the fact that the extended horseshoe doesn't has root, one self-contained proof for one-sided horseshoe  and another proof for two-sided horseshoe
  which is based on the no existence of root for the fullshift.

\medskip

\noindent {\bf Acknowledgement}
The authors would like to  thank   Klaus Schmidt, Evgeny Verbitskiy, Meng Wu and Weisheng Wu for helpful discussions and valuable comments. Thanks also go to Xiangdong Ye for providing the preprint \cite{HLXY2019} before it is diffused,  to B. Weiss for valuable informations.  Special thanks go to Evgeny Verbitskiy  for proposing to use an argument from Blum and Friedman \cite{BF1966} at the end of the proof of Proposition  6.1.
We are grateful to Shanghai Center for Mathematical Sciences for their hospitality, where part of this paper was written.

\section{Uniquely ergodic systems are not Bohr chaotic}
\label{Sect:UE}

We give here a quick proof of Theorem \ref{thm:UE}.
	It is a simple consequence of Robinson's topological weighted ergodic theorem, which states as follows (in a weak form)  \cite[Theorem 1.1]{Rob1994}.
Let $(X, T)$ be a uniquely ergodic dynamical system. If $\lambda\in \mathbb{S}^1 := \{z\in \mathbb{C}: |z| = 1\}$ is not an eigenvalue of the Koopman operator acting on $L^{2}(\mu)$ ($\mu$ being the invariant measure), then
for every continuous function $f\in C(X)$ we have
\begin{align*}
	\lim_{N\to \infty} \frac{1}{N} \sum_{k=0}^{N-1}\lambda^{-k}f(T^{k}x) =0,
\end{align*}
where the limit is actually uniform in  $x\in X$.
Note that eigenvalues are at most countable. Take
$\lambda\in \mathbb{S}^1$ which is not an eigenvalue.  The   nontrivial weight defined by  $w_{n}=\lambda^{-n}$ is then orthogonal
to the system $(X, T)$. Hence, the uniquely ergodic dynamical system $(X,T)$ is  not Bohr chaotic.
Theorem \ref{thm:UE} is thus proved.
\medskip

Let give some comments on systems which are not Bohr chaotic.

Krieger \cite{Krieger1972} proved that
every ergodic measure-preserving invertible  transformation of  a  Lebesgue  measure  space  is  isomorphic  to  an
orbital system on symbolic space, which is uniquely ergodic. So, by Theorem \ref{thm:UE}, every measure-preserving system has a topological model which is not  Bohr chaotic.
These models are Cantor sets.
\medskip

B\'eguin, Crovisier and Le Roux \cite{BCL2007} proved that
any compact manifold of dimension $d \ge 2$  which carries a minimal uniquely
ergodic homeomorphism $\mathcal{R}$ also carries a minimal uniquely ergodic homeomorphism $T$ with positive topological
entropy. Actually $T$ is an extension of $\mathcal{R}$.
By Theorem \ref{thm:UE}, such extensions are not Bohr chaotic.
 \medskip

  For any bounded, real-valued sequence $(w_n)$ with zero average along every infinite arithmetic progression (such sequences are said to be aperiodic), Downarowicz and Serafin \cite{DS2016} proved that there exists a subshift over $N$ symbols
to which $(w_n)$ is orthogonal and the entropy of the subshift
can approach to $\log N$. By Theorem 4 in \cite{Fan2017}, every $1$-oscillating sequence is aperiodic. Recall that $(w_n)$
is $1$-oscillating means
$$
  \forall t \in [0,1), \quad\lim_{N\to \infty}\frac{1}{N}\sum_{n=0}^{N-1} w_n e^{2\pi i n t}=0.
$$
That means, in harmonic analysis, that the Fourier-Bohr spectrum of $(w_n)$
is empty.
Note that for any real polynomial $P(x)=\sum_{k=0}^n a_k x^k\in \mathbb{R}[x]$
with at least one irrational coefficient   $a_k$ ($k\ge 2$),
$e^{2\pi i P(n)}$ is $1$-oscillating.

Also recall that Karagulyan \cite{Karagulyan2017} proved that
the M\"{o}bius sequence is not orthogonal to subshifts of finite type with positive topological entropy. But the Bohr-chaoticity requires more.
\medskip

\section{Subsystems of Bohr-chaotic systems}\label{Sect:subsystems}
We present here a basic idea for proving the Bohr chaoticity. It is to find a Bohr chaotic factor.

\subsection{Subsystems of Bohr-chaotic systems}
A system $(Y, S)$ is a {\em factor} of a system $(X,T)$ if there exists
a surjective continuous map $\pi: X \to Y$ such that $S\circ \pi = \pi \circ T$. In this case, $(X, T)$ is called an {\em extension} of $(Y, S)$. The following proposition follows from the definition.

\begin{prop} \label{factor} Any extension of a Bohr chaotic system is Bohr chaotic. Consequently, Bohr chaoticity is a topological invariant.
\end{prop}
\begin{proof}  Let $(Y,S)$ be a Bohr-chaotic factor of $(X, T)$
	with factor map $\pi$. Let $(w_n)$ be a non-trivial weight.
	Since $(Y,S)$ is Bohr-chaotic, there exists a
	pair $(g, y)$ with $g\in C(Y)$ and $y\in Y$ correlated to $(w_n)$
	for the system $(Y, S)$. Let
	$x$ be a pre-image of $y$ under $\pi$. Then $(g\circ \pi, x)$
	is a pair correlated to $(w_n)$ for the system $(X, T)$.  	
\end{proof}

The following proposition is also obvious.
\begin{prop} Let $\tau\ge 1$ be an integer. If $(X, T)$ is Bohr-chaotic, so is $(X, T^\tau)$.
\end{prop}

\begin{proof}
	Let $(w_n)$ be a non-trivial weight. Define a new one as follows: $v_{\tau n}=w_n$ and $v_j=0$ if $j$ is not a multiple of $\tau$. Suppose that $(f, x)$ is a pair correlated to $(v_j)$ for
	the system $(X, T)$. Then  $(f, x)$ is a pair correlated to $(w_n)$ for
	the system $(X, T^\tau)$.
\end{proof}

A partial inverse of the above proposition holds and it provides
a criterion for proving Bohr chaoticity.
Let $(X, T)$ be a given system and $K$ be a $T^\tau$-invariant compact set, i.e. $T^\tau K \subset K$. We say that $\tau$ is the {\em first return time} of $K$
if $T^k x \not\in K$ for all $x\in K$  and all $1\le k <\tau$.
Notice that $T^k K$ with $1\le k<\tau$ are all compact sets disjoint from $K$.
In this case, it is convenient to say that the $T^\tau$-invariant set $K$ has disjoint steps $T^jK$ ($1\le j<\tau$).
When $\tau=1$, it just means that $K$ is $T$-invariant.

\begin{thm}\label{ThmReturn} Let $(X, T)$ be a topological dynamical system. 
	Let  $K$ be a $T^\tau$-invariant compact subset ($\tau \ge 1$) having $\tau$ as its first return time.  If $(K, T^\tau)$ is Bohr-chaotic, so is $(X, T)$.
\end{thm}

\begin{proof} Let $(w_n)$ be a weight. Assume that
	$$
	\limsup_{N\to \infty}\frac{1}{N}\sum_{n=0}^{N-1}|w_{\tau n +j_0}|>0
	$$
	for some $0\le j_0 <\tau$. Since $(K, T^\tau)$ is Bohr chaotic,
	there exists $g \in C(K)$ and $x_0\in K$ such that
	\begin{equation}\label{sub1}
	\limsup_{N\to \infty}\frac{1}{N}\left|\sum_{n=1}^{N-1}w_{\tau (n-1) +j_0}g(T^{\tau n} x_0)\right|>0.
	\end{equation}
	Since $K^*:=\cup_{k=1}^{\tau-1}T^k K$ is disjoint from $K$, there
	exists a continuous function
	$g^*\in C(X)$ such that $g^*|_K=g$ and $g^*|_{K^*}=0$, by Urysohn's theorem. It follows that for $1\le j<\tau$, we have
	\begin{equation}\label{sub2}
	\lim_{N\to \infty}\frac{1}{N}\left|\sum_{n=1}^{N-1}w_{\tau (n-1) +j_0+j}g^*(T^{\tau n+j} x_0)\right|=0.
	\end{equation}
	From \eqref{sub1} and  \eqref{sub2} we get
	\begin{equation}\label{sub3}
	\limsup_{N\to \infty}\frac{1}{N}\left|\sum_{n=0}^{\tau N-1}w_{n +j_0}g^*(T^{n+j_0 +(\tau -j_0)} x_0)\right|>0.
	\end{equation}
	So, the function $f:=g^*$ and the point $x:=T^{\tau-j_0}x_0$
	satisfy the definition \eqref{BC} of the Bohr chaoticity
	of $(X, T)$.
\end{proof}

We shall prove Theorem \ref{thm:horse} by using Theorem \ref{ThmReturn} and the Bohr chaoticity of the full shift.

\subsection{Full shift}
On the symbolic space $\{0,1\}^\mathbb{N}$ we define the shift map
$\sigma$ by $(x_{n})_{n\ge 0} \mapsto (x_{n+1})_{n\ge 0}$.

\begin{prop}\label{FS} The one-sided full shift $(\{0,1\}^\mathbb{N}, \sigma)$ is Bohr chaotic.
\end{prop}

\begin{proof} Let $(w_n)$ be a non-trivial weight. We can assume that $w_n$'s are real numbers. Choose the function $f(x)= 1_{[0]}(x) - 1_{[1]}(x)$. Then choose the point $(x_n)$ defined by $x_n=0$ or $1$ according to $w_n\ge 0$ or $w_n <0$. Thus we have
\[
	   \limsup_{N\to \infty} \frac{1}{N}\sum_{n=0}^{N-1} w_n f(\sigma^n x) =
	   \limsup_{N\to \infty} \frac{1}{N}\sum_{n=0}^{N-1} |w_n| >0. \]
	\end{proof}

 The Bohr-chaoticity of the two-sided full shift  follows from that of the one-sided full shift.

\section{Affine maps on torus: Proof of Theorem \ref{torus}}
\label{Sect:Affine}

The proof of Theorem \ref{torus} was essentially contained in \cite{Fan1993}  in the context of uniform distribution on $\mathbb{T}^d$ of sequences of the form $(B^n x)$ where $B$ is an expanding matrix. It uses Riesz product measures on the group $\mathbb{T}^d$ to find points $x$ required in (\ref{BC})
and it is based on the lacunarity of the powers {  $B^n$. } We start with two combinatorial lemmas about the lacunarity.

\subsection{Two lemmas}
A sequence of vectors $H=(h_n)\subset \mathbb{R}^d$ ($d\ge 1$)
is said to be {\em dissociate} if any vector in $\mathbb{R}^d$
can be written in at most one way as a finite sum of the form
$$
\sum \epsilon_j h_j \quad {\rm with } \ \  \epsilon_j \in \{-1, 0, 1\}.
$$
This notion comes from Hewitt and Zuckermann \cite{HZ1966} and it allows us to define
the so-called Riesz product measures.
For any sequence $H$, we will use the following notation
$$
H^* =\left\{ \sum_{{\rm finite}}   \epsilon_j h_j:  \epsilon_j \in \{-1, 0, 1\}; h_j \in H \right\}.
$$

A sequence of complex numbers $\Lambda = (\lambda_n)_{n\ge 0}\subset \mathbb{C}\setminus \{0\}$
is said to be $\theta$-{\em lacunary} (\`a la Hadamard) for some  $\theta >1$ if
$|\lambda_{n+1}|\ge \theta |\lambda_n|$ for all $n \ge 0$.
A $\theta$-lacunary sequence with $\theta\ge 3$ is dissociate
(see a proof in \cite{Zygmund2002}, Vol. 1, p.208).

Given a lacunary sequence  of complex numbers $\Lambda=(\lambda_n)$ and an integer $q\ge 2$. We decompose $\Lambda$ into  $q$
parts in the following manner:
$$
\Lambda = \bigcup_{k=0}^{q-1}   \Lambda_k  \ \ {\rm with}\ \
\Lambda_k = (\lambda_{qn+k})_{n\ge 0}.
$$

Denote by $D(E, F)$ the distance between two sets $E$ and $F$.

\begin{lem}[\cite{Fan1993}] \label{lem:1} Suppose that $\Lambda=(\lambda_n)$ be a $\theta$-lacunary sequence of  complex numbers with $\theta >1$.
	Assume $0<\delta <|\lambda_0| (\theta-1)$. Then
	there exists an integer $q_0$ such that for $q\ge q_0$ we have\\
	\indent {\rm (i)}\  \  \ $\Lambda_0$ is dissociate;\\
	\indent {\rm (ii)}\  \
	$D(\Lambda_1 \cup \cdots \cup \Lambda_{q-1}, \Lambda_0^*) \ge \delta$;\\
	\indent {\rm (iii)} \ $D((\Lambda_p -\Lambda_q)\setminus \{0\}, \Lambda_0^*) \ge \delta$ for $1\le p \le q-1$.
	\\
	The above statements remain true if  $\Lambda_0$ is replaced by any
	other  $\Lambda_k$.
\end{lem}

Now let us consider a $d\times d$ matrix $B$ of entries in $\mathbb{Z}$ such that $\det B \not=0$. It defines an endomorphism
on $\mathbb{T}^d$. Let $B^*$ denote the transpose of $B$.
The matrix $B$ admits its complex Jordan normal form
$$
J =
\left[ {\begin{array}{cccc}
	J_1 &     &         &\\
	& J_2 &         &\\
	&     & \ddots &\\
	&     &        & J_r\\
	\end{array} } \right]
\ \ \ {\rm with} \ \ \
J_k = \left[ {\begin{array}{cccc}
	\xi_k &   1    &         &\\
	& \xi_k  & \ddots         &\\
	&        & \ddots & 1\\
	&     &        & \xi_k\\
	\end{array} } \right]
$$
where
$\xi_1, \xi_2, \cdots, \xi_r$ are the eigenvalues of $B$ with their multiplicities $d_1, d_2, \cdots, d_r$. Suppose that
$$
|\xi_1|\ge |\xi_2|\ge \cdots \ge |\xi_r|.
$$
We have
$$
B= T J T^{-1}
$$
where $T=(v_{1,1}, \cdots, v_{1,d_1}; \cdots; v_{r,1}, \cdots, v_{r,d_r})$ is a (complex) non singular matrix. The column vectors of $T$
are generalized eigenvectors of $B$:
$$
Bv_{k, 1}= \xi_k v_{k,1}, \quad
Bv_{k, j}= \xi_k v_{k,j} + v_{k, j-1} \
(1\le k\le r, 2\le j\le d_k).
$$

Choose a vector $h_0\in \mathbb{Z}^d$ such that $\langle v_{1,1}, h_0\rangle\not=0$. We can choose it among in the canonical basis
of $\mathbb{R}^d$. We are interested in the sequence
$H =(B^{*n} h_0)_{n\ge 0}$.  Let $h_n = B^{*n}h_0$.
Notice that
\begin{equation}\label{NN}
T^*\sum_{j=0}^n (\epsilon_j- \epsilon_j')h_j
= \sum_{j=0}^n (\epsilon_j- \epsilon_j')
{\rm diag} \{J_1^{*j}, \cdots, J_r^{*j}\} T^*h_0
\end{equation}
so that
\begin{equation}\label{NN2}
\|T^*\| \left\|\sum_{j=0}^n (\epsilon_j- \epsilon_j')h_j\right\|
\ge |\langle v_{1,1}, h_0\rangle|
\left|\sum_{j=0}^n (\epsilon_j- \epsilon_j')\xi_1^j
\right|.
\end{equation}
Here $T^*$ is considered as a linear operator on $\mathbb{C}^d$ which is equipped with the usual Hermite norm. To obtain the inequality
\eqref{NN2}, it suffice to estimate the first coordinate of the vector on the right hand side of \eqref{NN}.

Decompose $H$ into $q$ parts ($q\ge 2$):
$$
H = \bigcup_{k=0}^{q-1} H_k, \quad H_k =(h_{qn +k})_{n\ge 0}.
$$

From  the estimation \eqref{NN2} and Lemma \ref{lem:1}, we can get the following fact about the sequence of vectors $(h_n)$ which shares a refined dissociateness.

\begin{lem}[\cite{Fan1993}]\label{lem2} Let $\rho$ be the spectral radius of $B$
	and suppose $\rho>1$.
	For
	$0 <\delta <(\rho-1) |\langle v_{1,1}, h_0\rangle| /\|T\|$, there exists an integer $q_0\ge 2$ such that for $q\ge q_0$ we have\\
	\indent {\rm (i)} \ \ \ $H_0$ is dissociate;\\
	\indent {\rm (ii)} \  \ $D(H_1 \cup \cdots \cup H_{q-1}, H_0^*) \ge \delta$;\\
	\indent {\rm (iii)} \ $D((H_p -H_q)\setminus \{0\}, H_0^*) \ge \delta$ for $1\le p \le q-1$.
	\\
	The above statements remain true if  $H_0$ is replaced by any $H_k$.
\end{lem}

\subsection{Proof of Theorem \ref{torus}}
For $Tx = Bx +b$ we have
$$
T^n x = B^n x + (B^{n-1}+ \cdots+B +I)b.
$$
Take $f(x) = e^{2\pi i \langle h_0, x\rangle}$. We have
$$
f(T^n x) = e^{2\pi i \psi_n} e^{2\pi i \langle B^{*n}h_0, x\rangle}
$$
where
$ \psi_n = \langle h_0, (B^{n-1}+ \cdots+B +I)b \rangle
$.

Let $q$ be sufficiently large such that the statements in Lemma
\ref{lem2} hold. For any given non-trivial weight $(w_n)$, there
exits $0\le k\le q-1$ such that
\begin{equation}
\limsup_{N\to \infty}\frac{1}{N}\sum_{n=0}^{N-1}|w_{qn+k}|>0.
\end{equation}
Without loss of generality, we assume that $k=0$. Since
$H_0$ is dissociate (Lemma \ref{lem2} (i)), for any sequence
of complex numbers $(a_n)$ such that $|a_n|\le 1$, we can define
the following Riesz product measure on $\mathbb{T}^d$ (see \cite{HZ1966})
\begin{equation}
\nu_a = \prod_{n=0}^\infty \left(1 + {\rm Re} (a_n\ e^{2\pi i \langle B^{*qn} h_0, x\rangle})\right).
\end{equation}
This Borel probability measure $\nu_a$ is characterized by its Fourier
coefficients as follows
\begin{equation}\label{F1}
\widehat{\nu}_a(\sum \epsilon_n B^{*qn}h_0) = \prod a_n^{(\epsilon_n)}
\end{equation}
where $a_n^{(0)} =1, a_n^{(1)} = a_n/2$ and $a_n^{(-1)} = \overline{a}_n/2$; and
\begin{equation}\label{F2}
\widehat{\nu}_a({\rm \textbf{n}})=0
\quad {\rm if}\ \  {\rm \textbf{n}} \not \in H_0^*.
\end{equation}

By the formula \eqref{F1}, $\{e^{2\pi i \langle B^{*qn}h_0, x\rangle} -\frac{\overline{a}_n}{2}\}$ is an orthogonal system in $L^2(\nu_a)$. Then, by the Menshov-Rademacher theorem (see \cite{Zygmund2002} vol. p. 193. see also \cite{Fan1993b} ),
the series $$
\sum \frac{w_{qn}e^{2\pi i \psi_{qn}}}{n}\left(e^{2\pi i \langle B^{*qn}h_0, x\rangle} -\frac{\overline{a}_n}{2}\right)$$
converges $\nu_a$-almost everywhere. By Kronecker's lemma, it follows that
$\nu_a$-almost everywhere
\begin{equation*}
	\lim_{N\to \infty}\frac{1}{N}\sum_{n=0}^{N-1} w_{qn} e^{2\pi i \psi_{qn}}
	\left(e^{2\pi i \langle B^{*qn}h_0, x\rangle} -\frac{\overline{a}_n}{2}\right) =0.
\end{equation*}

We will fix our choice for $a_n$ as follows
\begin{equation}\label{an}
a_n = r e^{i\arg w_{qn} + 2\pi i \psi_{qn}},
\end{equation}
where $r$ is any fixed number such that $0<r\le 1$. So, for $\nu_a$-almost every $x$ we have
\begin{equation}\label{C1}
\limsup_{N\to \infty}\frac{1}{N}\left|\sum_{n=0}^{N-1} w_{qn} f(T^{qn} x)\right| =
\frac{r}{2}\limsup_{N\to \infty} \frac{1}{N} \sum_{n=0}^{N-1} |w_{qn}|
>0.
\end{equation}

Fix  $1\le p<q$ and let $X_n = e^{2\pi i \langle B^{*qn+p}h_0, x\rangle}$.
By the formula \eqref{F2} and Lemma \ref{lem2} (ii) and (iii),  we have
$$
\mathbb{E}_{\nu_a} X_n =0, \quad  \mathbb{E}_{\nu_a} X_n \overline{X}_m =0 \ \ {\rm if} \ \  n\not=m.
$$
Again, by using the Menshov-Rademacher theorem to the orthogonal
system $\{X_n\}$, for $\nu_a$-almost every $x$ we have
\begin{equation}\label{C2}
\lim_{N\to \infty}\frac{1}{N}\sum_{n=0}^{N-1} w_{qn+p} f(T^{qn+p} x) = 0.
\end{equation}

From \eqref{C1} and \eqref{C2} we conclude that
for $\nu_a$-almost every point $x$ we have
\begin{equation} \label{C3}
\limsup_{N\to \infty}\frac{1}{N}\left|\sum_{n=0}^{N-1} w_{n} f(T^{n} x)\right| > 0.
\end{equation}
Thus we have proved the Bohr chaoticity of  $(\mathbb{T}^d, T)$.

The Hausdorff dimension of the set of $x$ such that
\eqref{C3} holds is not less than  the Hausdorff dimension of the measure
$\nu_a$. But $\dim \nu_a$ tends to $d$ when $r \to 0$. 	
\qed 	
\bigskip

We remark that if $(w_n)$ is aperiodic and if we take $a_n = r e^{2\pi i \psi_n}$, the
weighted ergodic limit exists:
\begin{equation*}
	\nu_a\!-\!a.e. \quad\lim_{N\to \infty}\frac{1}{N}\sum_{n=0}^{N-1} w_{n} f(T^{n} x) = \frac{r}{2} \lim_{N\to \infty}\frac{1}{N}\sum_{n=0}^{N-1} w_{qn}.
\end{equation*}
Notice that the above limit on the right hand side exists since $(w_n)$
is aperiodic.

\subsection{Bohr chaotic automorphisms on the torus having no horseshoe}
Hyperbolic automorphisms on torus admit horseshoes. So, Theorem \ref{thm:horse} covers the class of all hyperbolic automorphisms on torus. But it doesn't cover all endomorphisms of positive entropy concerned by Theorem \ref{torus}. For example, let $B\in {\rm GL}(d, \mathbb{Z})$ with $\det B = \pm 1$. Suppose that \\
\indent {\rm (i)} \ the characteristic polynomial $\chi_B$ of $B$ is irreducible over $\mathbb{Q}$;\\
\indent {\rm (ii)} some but not all of eigenvalues of $\chi_B$ is on the unit circle. \\
Then the automorphism on the torus $\mathbb{T}^d$ defined by $Tx = Bx \mod \mathbb{Z}^d$ has no homoclinic point (see Example 3.4. in \cite{LS1999}).
Consequently, for any integer $m\ge 1$, the power $T^m$ has no homoclinic point.
Since the full shift $(\{0,1\}^\mathbb{Z}, \sigma)$ has a dense set of homoclinic points,  $T$ can not admit a horseshoe.

Therefore it is not possible to prove Theorem \ref{torus} by constructing a horseshoe as we do for proving Theorem \ref{thm:horse}. The family of systems admitting horseshoes is a proper sub-family of the family of Bohr chaotic systems.

However, the above automorphism $T$ admits a semi-horseshoe, because according to \cite{HLXY2019}, any automorphism on a compact metric abelian group  having positive
topological entropy  has a semi-horseshoe. The above  automorphism also shows that a partially hyperbolic system can have a semi-horseshoe, but not a horseshoe.

\section{Proof of Theorem \ref{thm:horse}: case of one-sided horseshoe}\label{sect:1horseshoe}
In order to prove Theorem \ref{thm:horse} (the case of one-sided horseshoe),  it suffices to prove,  by Theorem \ref{ThmReturn},  that
a topological dynamical system having a one-sided horseshoe  must have a horseshoe with disjoint steps, as we state as follows

\begin{thm} \label{thm:disjiontHS}
	Suppose that $(X,T)$ is a topological dynamical system admitting a one-sided $N$-order horseshoe $(\Lambda, T^N)$.
	Then there exist a  closed subset $K\subset \Lambda$ and a positive integer $M$ such that $(K,T^{MN})$ is a  one-sided horseshoe of $(X,T)$  with disjoint steps.
\end{thm}

The strategy for proving Theorem \ref{thm:horse} in the case of two-sided horseshoe is the same.   But there some differences in details.
The proof of this case is postponed to the next section.

The proof of Theorem \ref{thm:disjiontHS} is rather long. In the next, we first present the ideas of proof and then present the proof in several subsections.

\subsection{Ideas of proof}
 Let $(\Lambda, T^N)$ be the given horseshoe.  We consider the extended horseshoe
$T: \Lambda^*\to \Lambda^*$ where
$$
    \Lambda^* = \bigcup_{n=0}^{N-1} T^n(\Lambda).
$$
It is a subsystem of $(X, T)$.
One of key points for proving Theorem \ref{thm:disjiontHS} is that
$(\Lambda^*, T)$ has no roots, i.e. there is no continuous map $S:  \Lambda^*\to \Lambda^*$ such that $S^n=T$ for some $n\ge 2$.
Its exact statement is as follows.

\begin{prop}\label{prop:no-root} Suppose that $(X,T)$ is a topological dynamical system admitting a $N$-order one-sided horseshoes $(\Lambda, T^{N})$ such that\\
	\indent {\rm (i)} $(\Lambda, T^{N})$ has disjoint steps;\\
	\indent {\rm (ii)} every map $T:T^{j}(\Lambda)\to T^{j+1}(\Lambda)$ is  bijective for  $0\leq j\leq N-2$. \\
	Then for any integer  $n\geq 2$, there is no continuous map   $S:\Lambda^*\to \Lambda^*$ such that $S^n=T$.
\end{prop}
The proof of this proposition will be given in Appendix A.

The second key point  is
that in any cylinder of $\{0, 1\}^{\mathbb{N}}$,  
there is a horseshoe with disjoint steps for the shift map. This will be proved in Appendix B.
Then,
since $(\Lambda, T^N)$ is conjugate to $\{0, 1\}^{\mathbb{N}}$, we can find
 horseshoes $\Lambda'$ in any non-empty open set of  $\Lambda$. These horseshoes $\Lambda'$, called prototypes of horseshoes,   have not yet the required disjoint steps, but some steps are really disjoint from $\Lambda'$.  Actually, from the fact that $(\Lambda, T^N)$ is conjugate to the shift map, we can easily find our
first prototype of horseshoe which is of the form $(\Lambda', T^{NM})$ for some integer $M\ge 1$ such that $T^{kN}(\Lambda')$ for $1\le k <M$ are disjoint from
$\Lambda'$ (a partial disjointness). This is the starting point of our construction by induction.

In general, given a horseshoe $(\Lambda', T^{pq})$ ($p, q$ being integers) such that
\begin{equation}\label{eq:Disj1}
    \Lambda' \cap T^{kq}(\Lambda')= \emptyset \qquad (1\le k <p),
\end{equation}
we have
\begin{equation}\label{eq:Disj2}
     \Lambda' \not\subset \bigcup_{k=0}^{p-1} T^{kq +s}(\Lambda') \text{ for any $1\le s<q$ with $s|q$ }
\end{equation}
(cf. Proposition \ref{PO2}). This is a consequence of the fact that $T$ has no roots. We could say that $\Lambda'$ is displaced by
$T^{kq +s}$'s.  The second key point allows us to find
a better horseshoe $\Lambda'' \subset \Lambda'$ such that
\begin{equation}\label{eq:Disj3}
\Lambda'' \subset   \Lambda'\setminus \bigcup_{k=0}^{p-1} T^{kq +s}(\Lambda').
\end{equation}
Here, by ``better'' we mean that $\Lambda''$ has ``more'' disjointness. Indeed, from $(\ref{eq:Disj1})$, $(\ref{eq:Disj2})$  and $(\ref{eq:Disj3})$ we get
$$
\Lambda'' \cap T^{kq}(\Lambda'')= \emptyset,
\quad
\Lambda'' \cap T^{kq+s}(\Lambda'')= \emptyset \qquad (1\le k <p,
1\le s <q \ {\rm with}\  s|q).
$$

Basing on the decomposition of $N=p_1^{\gamma_1}\cdots p_s^{\gamma_s}$
into primes,  recursively we can find finer and finer horseshoes  to finally get what we want.

We invite the readers to read the isolated proofs for the case $N=2$ (Subsection \ref{sec:N=2}) and the case  $N=3$  (Subsection \ref{sec:N=3}), which show some more details of the above ideas.   The machinery of induction is explained
in Proposition \ref{prop:disp1} and Proposition \ref{prop:disp2}.
The proof in the case $N=2$ doesn't need Proposition \ref{prop:disp1} and Proposition \ref{prop:disp2}. The proof in the case $N=3$ only need Proposition \ref{prop:disp1}.

\subsection{Displacement of horseshoe}

Let   $(\Lambda,T^N)$ be a $N$-order one-sided  horseshoe.
If $\Lambda \not\subset T^j(\Lambda)$ for some $1\le j <N$, we say that
$\Lambda$ is displaced by $T^j$. We shall prove that the horseshoe $\Lambda$ is displaced by $T^j$ when $j|N$. Actually the following more general  Proposition \ref{PO2} , which
is a consequence of Proposition \ref{prop:no-root}. 

\begin{prop}\label{PO2}
	Let  $(X, T )$  be  a topological dynamical system. Suppose that  $(X, T)$ has a  $pq$-order one-sided horseshoe $(\Lambda, T^{pq} )$ for some integers $p\ge 1$ and $q\ge 1$ such that \\
	\indent {\rm  (i)} \ $\Lambda\cap T^{kq}(\Lambda)=\emptyset$ for $1\leq k\leq p-1$;\\
	\indent {\rm  (ii) } $T^{(p-1)q}:\Lambda\to T^{(p-1)q}(\Lambda)$ is a bijection. \\
	Then for any integer $1\le s < q$ such that  $s\mid q$ we have $ \Lambda\not\subset \bigcup_{k=0}^{p-1}T^{kq+s}(\Lambda)$.
\end{prop}
\begin{proof}
	We give a proof by contradiction.  Suppose  that $\Lambda \subset \bigcup_{k=0}^{p-1}T^{kq+s}(\Lambda)$ for some
	$s$ with $s |q$. Then, from the fact $\Lambda =T^{pq}(\Lambda)$, we get
	\[T^{\ell q}(\Lambda)\subset T^{\ell q} \left(\bigcup_{k=0}^{p-1}T^{kq+s}(\Lambda)\right)=\bigcup_{k=0}^{p-1}T^{kq+s}(\Lambda) \quad \text{for } 1\leq \ell \leq {p-1}.\] Therefore,  $\bigcup_{k=0}^{p-1}T^{kq}(\Lambda) \subset \bigcup_{k=0}^{p-1}T^{kq+s}(\Lambda)$.
	Let $$\Lambda_q^*=\bigcup_{k=0}^{p-1}T^{kq}(\Lambda).
	$$ We claim that $\Lambda_q^* = T^s (\Lambda_q^*)$. Indeed, from $\Lambda_q^* \subset T^s (\Lambda_q^*)$ which is just proved above,  we have
	\[\Lambda_q^* \subset T^s \big(\Lambda_q^*\big)\subset T^{2s} \big(\Lambda_q^*)\subset \cdots \subset
	T^{pq}(\Lambda_q^*)=\Lambda_q^*.\]

	Now let us consider the  dynamical  system $(\Lambda_q^*, T^q)$ which admits a $p$-order horseshoe $(\Lambda, T^{pq})$ with disjoint steps.
	On one hand, $T^s:\Lambda_q^*\to \Lambda_q^*$ is a $q/s$-root of $T^{q}$. On the other hand, the following factor maps are   bijective:
	\begin{equation}\label{eq:factormaps}
	\Lambda \stackrel{T^q}{\longrightarrow}  T^{q}(\Lambda)
	\stackrel{T^q}{\longrightarrow}  T^{2q}(\Lambda) \stackrel{T^q}{\longrightarrow} \cdots
	\stackrel{T^q}{\longrightarrow}  T^{(p-2)q}(\Lambda)  \stackrel{T^q}{\longrightarrow}  T^{(p-1)q}(\Lambda)
	\end{equation}
	because  the map $T^{(p-1)q}:\Lambda\to T^{(p-1)q}(\Lambda)$ is supposed bijective. These two facts together
	contradict to  Proposition  \ref{prop:no-root} when we consider the system $(\Lambda_q^*, T^q)$.
	\end{proof}

Taking $p=1$ and $N=q$, as an immediate consequence, we obtain that $\Lambda\not \subset T(\Lambda)$.
\begin{cor}\label{PO1}
	Let $(X, T )$ be a topological dynamical system admitting an  $N$-order horseshoe $(\Lambda, T^N )$ with $N\geq 2$. Then for any positive integer $1\le s< N$ dividing $N$ we have $\Lambda\not \subset  T^s(\Lambda)$.
\end{cor}

\medskip

\subsection{Proof of Theorem \ref{thm:disjiontHS}: $N=2$}
\label{sec:N=2}

The proof will be based on
Proposition \ref{PO2} (displacement of horseshoe) and   Proposition \ref{lem-SBSYS} (construction of horseshoe in a cylinder). The argument will be inductive and several constructions of horseshoe will be performed and the final horseshoe will have a high order. The proof will be relatively easy if the order $N$ is a power of  prime i.e. $N=p^m$. In the following, we first give a proof for $N=2$ which is rather direct and for $N=3$ which shows the idea for the general case.
\medskip

{\it Proof for $2$-order horseshoe.} Suppose $(\Lambda, T^2)$ is an one-side  horseshoe of the topological system $(X, T)$.  By Corollary \ref{PO1}, $\Lambda$ is displaced by $T$. In other words, there exist a point $x\in \Lambda$  such that $x\notin T(\Lambda)$.   Since $T$ is continuous and  $\Lambda$ is compact, there exist a neighborhood $U\subset \Lambda$ of $x$ such that $U\cap T (\Lambda)=\emptyset$.
Since $(\Lambda, T^2)$ is topologically conjugate to the full shift system $(\{0,1\}^{\N},\sigma)$, by Proposition \ref{lem-SBSYS}, there exist  an integer $M$ and a  subset $K\subset U$ such that $(K,T^{2M})$ is a horseshoe  with
$$
K\cap T^{2j}(K)=\emptyset \ \ \ {\rm for}\ \  1\leq j\leq M-1.
$$
Notice  that $T^{2j+1}(K)\subset T^{2j+1}(\Lambda)=T(\Lambda)$,
because of $T^2(\Lambda)=\Lambda$. From $U\cap T(\Lambda)=\emptyset$, it follows that
$$
K\cap T^{2j+1}(K)=\emptyset \ \ \ {\rm for}  \ \ j\ge 0.$$
Thus we have proved that $(K,T^{2M})$ is a horseshoe  with disjoint steps of $(X,T)$.

\subsection{Proof of Theorem \ref{thm:disjiontHS}: $N=3$}
\label{sec:N=3}
{\it Proof for $3$-order horseshoe.}  Suppose $(\Lambda, T^3)$ is a one-side  horseshoe of the topological system $(X, T)$. The proof is decomposed into  two steps.
\begin{itemize}
	\item By the same argument as in the case for $2$-order horseshoe, there exist an integer $M_1$ and  a  subset $\Lambda_1$ such that  $(\Lambda_1,T^{3M_1})$ is a horseshoe with
	\[
	\Lambda_1\cap T^{3j}(\Lambda_1)=\emptyset  \quad \text{for } 1\leq j\leq M_1-1,
	\]
	and
	\[
	\Lambda_1\cap T^{3j+1}(\Lambda_1)=\emptyset  \quad \text{for } 0\leq j\leq M_1-1.
	\]
	\item This second step is sketchy. Details are in the proof of Proposition  \ref{prop:disp1} (corresponding to the special case $p=M_1,q=3, J=\{1\}$ ).  In general, $\Lambda_1$ may intersect $T^{3j+2}(\Lambda_1)$ for some $0\leq j \leq M_1-1$.
	We claim the following displacement (see (\ref{eq:CNC1})):
	$$
	\Lambda_1\not\subset \bigcup_{j=0}^{M_1-1}T^{3j+2}(\Lambda_1).
	$$ Then we can construct a
	subset $\Lambda_2\subset \Lambda_1$ such that $(\Lambda,T^{3M_1M_2})$ is a horseshoe for some  integer $M_2$. It can be checked that the horseshoe   $(\Lambda,T^{3M_1M_2})$ has disjoint steps.
\end{itemize}

{Here, for the second step, we prove the displacement  by  using  the  simple arithmetic fact $2+2 \equiv 1 (\!\!\!\!\mod 3)$(taking $s=2$, $n=2$ and $q=3$ in Proposition \ref{prop:disp1}).}
\medskip

\subsection{Construction of finer horseshoes}
The following proposition improves the disjointness of the steps of a horseshoe, by construction a smaller horseshoe.

\begin{prop}\label{prop:disp1}
	Let  $p,q $ be  two positive  integers and   let  $J\subset \{1, 2, \cdots, q\}$. Suppose that a topological system $(X, T)$ has a $T^{pq}$-invariant subset $\Lambda\subset X$  such that \\
	\indent {\rm  (a)} $(\Lambda, T^{pq})$ is a horseshoe;\\
	\indent {\rm  (b)} $\Lambda\cap T^{kq}(\Lambda)=\emptyset$ for $1\leq k\leq p-1$;\\
	\indent {\rm  (c)} $\Lambda\cap T^{kq+j}(\Lambda)=\emptyset$ for $0\leq k\leq p-1$ and $j\in J$\\
	Then for any integer $s\ge 1$ such that  $ns\ (\!\!\!\!\mod q) \in J$ for some positive integer $n$, there exist an integer $M\ge 1$ and a subset $\Lambda^{\prime}\subset \Lambda$
	such that \\
	\indent {\rm  (A)} $(\Lambda^{\prime},T^{Mpq})$ is a horseshoe; \\
	\indent {\rm  (B)} $\Lambda^{\prime}\cap T^{kq}(\Lambda^{\prime})=\emptyset$ for $1\leq k\leq pM-1$;\\
	\indent {\rm  (C)} $\Lambda^{\prime}\cap T^{kq+j}(\Lambda^\prime)=\emptyset$ for $0\leq k\leq pM-1$ and $j\in J\cup \{s\}$.
\end{prop}
\begin{proof}
	Let  $\Lambda^*=\bigcup_{k=0}^{p-1}T^{kq}(\Lambda)$.
	We claim that
	\begin{equation}\label{eq:CNC1}
	\Lambda\not\subset  T^{s}\left(\Lambda^*\right).
	\end{equation}
	Otherwise, by the fact $T^{pq}(\Lambda)=\Lambda $, we would  have $T^{q}(\Lambda^*)=\Lambda^*$, which implies \[T^{kq}(\Lambda)\subset T^{kq+s}(\Lambda^*)=T^{s}(\Lambda^*)\text{ for } 0\leq k\leq  p-1. \]
Hence, $\Lambda^*\subset  T^{s}(\Lambda^*)$,  and consequently
	\begin{equation}\label{eq:E11}
	\forall \ell \in \mathbb{N},\quad \Lambda^*\subset  T^{\ell s}\left(\Lambda^*\right).
	\end{equation}
	However, by the assumption on $s$, there exists a positive  integer $n$ such that
	$$
	T^{ns}\left(\Lambda^*\right)\subset \bigcup_{j\in J} T^{j}\left(\Lambda^*\right).$$
	This, together with (\ref{eq:E11}), implies
	\[\Lambda^*\subset  \bigcup _{j\in J} T^{j}\left(\Lambda^*\right).\]
	This contradicts  the condition (c).
	
	By (\ref{eq:CNC1}), there exists a point  $x\in \Lambda \setminus T^{s}\left(\Lambda^*\right)$. We  take an open neighborhood $U$ of $x$
	such that $U\cap  T^{s}(\Lambda^*)=\emptyset$.
	By  the condition  (a) and Proposition \ref{lem-SBSYS},
	there exists  a $T^{Mpq}$-invariant set $\Lambda^{\prime}\subset \Lambda\cap  U$ for some integer $M\ge 1$   such that
	$(\Lambda^{\prime},T^{Mpq})$ is conjugated to $(\{0,1\}^\N,\sigma)$ and
	\begin{align}\label{eq:s}
		\Lambda^{\prime}\cap T^{\ell pq}(\Lambda^\prime)=\emptyset \quad \text{ for }1\leq \ell \leq M-1.
	\end{align}

	It remains to check (B) and (C).
	
	Write
	$k=sp+r$ with $1\le s \le M-1$ and $0\le r\le p-1$.
	If $r=0$, we get the property (B) from (\ref{eq:s}).
	If $r\not=0$,  we get the property (B) from (b) and the facts
	$\Lambda^\prime \subset \Lambda$ and  $ T^{pq}\Lambda=\Lambda$. Indeed,
	\[
	\Lambda^\prime \cap T^{kq}(\Lambda^\prime)=\Lambda^\prime \cap T^{s pq + r q}(\Lambda^\prime)\subset \Lambda^\prime \cap T^{rq}(\Lambda)  =\emptyset.
	\]
	(B) is thus checked.
	
	Now we first check (C) for $j \in J$, using (c).
	Write
	$k=sp+r$ with $0\le s \le M-1$ and $0\le r\le p-1$.
	Using $T^{pq}(\Lambda) =\Lambda$ and $\Lambda'\subset \Lambda$, we get
	$$
	T^{kq+j}(\Lambda') \subset T^{spq + rq +j}(\Lambda)
	= T^{rq +j}(\Lambda).
	$$
	The last set is disjoint from $\Lambda$, by (c). So, it is disjoint from $\Lambda'$.
	Then let us check (C) for $j =s$.
	Note that
	$$
	T^{k q +s}(\Lambda') \subset  T^{k q +s}(\Lambda)
	\subset T^s(\Lambda^*).
	$$
	The facts $\Lambda^{\prime} \subset U$ and $U \cap  T^{s}(\Lambda^*)=\emptyset$ imply that
	$\Lambda'$ is disjoint from $\bigcup_{k=0}^{p-1} T^{kq +s}(\Lambda)$. This, together with $T^{pq}(\Lambda)=\Lambda$,
	implies that $\Lambda'$ is disjoint from $\bigcup_{k=0}^{Mp-1} T^{kq +s}(\Lambda^\prime)$.
\end{proof}

The proof for the case of  $3$-order horseshoe uses  Proposition \ref{prop:disp1}. Using Proposition \ref{prop:disp1} inductively,  we can give a complete proof for the case  $N=p^\gamma$.    

However, for the general case, Proposition \ref{prop:disp1} is not enough for the proof of Theorem \ref{thm:disjiontHS}.
For example, when $N=6$, for any factor  $j$ of $6$, there exists a number $s\in \{2,3\}$ such that
	\[ns\not\equiv j (\!\!\!\!\mod 6) \text{ for any positive integer  $n$}.  \] 
The following proposition is another improvement of disjointness of steps of horseshoe.

\begin{prop}\label{prop:disp2}
	Let  $p,q $ be  two positive  integers and   let  $J\subset \{1, 2, \cdots, q\}$. Suppose that a topological system $(X, T)$ has a $T^{pq}$-invariant subset $\Lambda\subset X$  such that  \\
	\indent {\rm  (a)} $(\Lambda, T^{pq})$ is a horseshoe;\\
	\indent {\rm  (b)} $T^{(p-1)q}:\Lambda\to T^{(p-1)q}(\Lambda)$ is a bijection; \\
	\indent {\rm  (c)} $\Lambda\cap T^{kq}(\Lambda)=\emptyset$ for $1\leq k\leq p-1$;\\
	\indent {\rm  (d)} $\Lambda\cap T^{kq+j}(\Lambda)=\emptyset$ for $0\leq k\leq p-1$ and $j\in J$. \\
	Then for any integer $s\ge 1$ such that  $s\mid q$, there exist an integer $M$ and a subset $\Lambda^{\prime}\subset \Lambda$
	such that \\
	\indent {\rm  (A)} $(\Lambda^{\prime},T^{Mpq})$ is a horseshoe; \\
	\indent {\rm  (B)} $T^{(pM-1)q}:\Lambda^\prime\to T^{(pM-1)q}(\Lambda^\prime)$ is a bijection;\\
	\indent {\rm  (C)} $\Lambda^{\prime}\cap T^{kq}(\Lambda^{\prime})=\emptyset$ for $1\leq k\leq pM-1$;\\
	\indent {\rm  (D)} $\Lambda^{\prime}\cap T^{kq+j}(\Lambda^\prime)=\emptyset$ for $0\leq k\leq pM-1$ and $j\in J\cup \{s\}$.
\end{prop}

\begin{proof}Let  $\Lambda^*= \bigcup_{k=0}^{p-1}T^{kq}(\Lambda)$.
	We have $T^{s}(\Lambda^*)=\bigcup_{k=0}^{p-1}T^{kq+s}(\Lambda)$.
	Under the condition (a), (b) and (c), we apply Proposition  \ref{PO2} to get the displacement  $ \Lambda\not\subset T^{s}(\Lambda^*)$.  Choose  $x\in \Lambda$ such that $x\notin T^{s}(\Lambda^*)$ and then choose
	an open neighborhood $U$ of $x$
	such that $U$ and $T^{s}(\Lambda^*)$ are disjoint.  As $(\Lambda, T^{pq} )$ is conjugate to the full shift $(\{0,1\}^{\N},\sigma)$, we can assume that $U=\pi (C)$ for some  cylinder $C$, where   $\pi: (\{0,1\}^{\N},\sigma) \to
	(\Lambda, T^{pq} ) $ is the conjugation.  Applying  Proposition \ref{lem-SBSYS}
	to $C$ to get a horseshoe contained in $C$ and then projecting  the horseshoe by $\pi$. Thus we get a $T^{Mpq}$-invariant set $\Lambda^{\prime}\subset U\cap \Lambda$ for some positive integer $M$,
	such that
	\begin{itemize}
		\item[(A)] $(\Lambda^{\prime},T^{Mpq})$ is a horseshoe;
		\item[(B')] $T^{(M-1)pq}: \Lambda'\to  T^{(M-1)pq}(\Lambda')$ is bijective;
		\item[(C')] $\Lambda^{\prime}\cap T^{kpq}(\Lambda')=\emptyset$ for $1\leq k\leq M-1$;
		\item[(D')] $\Lambda^{\prime}\cap T^{kq+s}(\Lambda^\prime)=\emptyset$ for  $0\leq k\leq pM-1$.
	\end{itemize}
	Let us check (C),  (D) and then (B).
	Since $\Lambda'\subset \Lambda$ and $T^{pq}(\Lambda)=\Lambda$, (C) follows from  (c) and (C').
	Since  \[\Lambda^{\prime} \cap  T^{s}(\Lambda^*)=\emptyset,\] (D) follows from (d).
	Recall the facts $\Lambda^{\prime}\subset \Lambda$ and $T^{(M-1)pq} (\Lambda)=\Lambda$.  Also recall that $T^{(p-1)q}$ is a bijection from $\Lambda$ to $T^{(p-1)q}(\Lambda)$, by (b). It follows that   $T^{(p-1)q}$ is a bijection from $T^{(M-1)pq}(\Lambda^{\prime})$ to $T^{(Mp-1)q}(\Lambda^{\prime})$. Hence, (B) follows  from (B').
\end{proof}

\subsection{A simple arithmetic fact}
In the proof of Theorem \ref{thm:disjiontHS}, we shall use a simple arithmetic fact stated as Lemma \ref{lem-final} below.   Let $N\ge 2$ be an integer and let us consider the cyclic group   $\mathbb{Z}/N\mathbb{Z}$ identified with  $\{0,1,\cdots N-1\}$. For each element $a\in \mathbb{Z}/N\mathbb{Z}$,  $\langle a\rangle $ will denote the subgroup generated by $a$. Assume that $N=p_{1}^{\gamma_1}p_{2}^{\gamma_2}\cdots p_{s}^{\gamma_s}$ is  decomposed into primes with $\gamma_j\ge 1$ ($j=1,2,\cdots s$).
For $n\in \Z/N\Z$,  denote by $\langle n \rangle$  the subgroup of $\Z/N\Z$ generated by $n$. It is easy to see that the subgroups  $\langle N/p_1\rangle, \cdots, \langle N/p_s\rangle$ are exactly the  smallest nontrivial subgroups  of $\mathbb{Z}/N\mathbb{Z}$. Here, by smallest nontrivial subgroup, we mean a subgroup which  does not contain a non-trivial proper subgroup.   It is also easy to see that each nontrivial subgroup of $\mathbb{Z}/N\mathbb{Z}$ contains at least one subgroup  $\langle N/p_i\rangle$.
These facts imply the following lemma.

\begin{lem}\label{lem-final}
	
	Let  $N=p_{1}^{\gamma_1}p_{2}^{\gamma_2}\cdots p_{s}^{\gamma_s}$ be a positive integer  decomposed into primes with $\gamma_j\ge 1$ ($j=1,2,\cdots s$).
	Then for each integer $n\in \{1,2,\cdots, N-1\}$, there exist an integer $x \in \{1,2,\cdots, N-1\}$ and a prime $p_i$ ($1\le i \le s$) such that
	$$n x = \frac{N}{p_i} \mod N.$$ \end{lem}

\begin{proof}
	For each given nonzero element $n\in \mathbb{Z}/N\mathbb{Z}$, the group $\langle n\rangle$ contains  a subgroup  $\langle N/p_i\rangle$ for some $1\leq i \leq s$. So,   $N/p_i \in \langle n\rangle $, which means that
	$ N/p_i = kn \mod  N$
	for some  positive integer $k\in\{1,2,\cdots, N-1\}$.
\end{proof}

\subsection{Proof of Theorem \ref{thm:disjiontHS}: for general $N\ge 1$}

There is nothing to prove if $N=1$. Then assume $N\ge 2$.
	Write $N=p_{1}^{\gamma_1}p_{2}^{\gamma_2}\cdots p_{s}^{\gamma_s}$, where $p_1>p_2>\dots >p_s\ge 2$ are distinct prime numbers and $\gamma_1  \ge 1, \cdots, \gamma_s \ge 1$.
	For $1\leq i \leq s$, let  $N_i=\frac{N}{p_i}$.

	We shall divide our proof into two steps. Our machinery is to construct a sub-horseshoe whenever a horseshoe is displaced, using Proposition \ref{prop:disp1} or Proposition \ref{prop:disp2}, which are based on  Proposition \ref{PO2} and Proposition \ref{lem-SBSYS} in Appendix B. The fist step is to construct  a horseshoe  $(\Lambda_s, T^{L_sN})$, where $\Lambda_s \subset \Lambda$ and  $L_s > 1$ is an integer,
	such that \\
	\indent (1) \ $\Lambda_s\cap T^{nN}(\Lambda_s)=\emptyset$  for   $1\leq n \leq L_s-1$;\\
	\indent (2) \ $\Lambda_s\cap T^{nN+N_i}(\Lambda_s)=\emptyset$    for all $1\leq i\leq s$ and
	all $0\leq n\leq L_s-1$.  \\
	The second step is to repeat the same machinery to construct a subset $ \Lambda^{\prime} \subset \Lambda_s$ such that  for some integer $M$ we have \\
	\indent (3)\    $(\Lambda^{\prime}, T^{MN})$ is a  horseshoe; \\
	\indent (4) \
	\( \Lambda^{\prime}\cap T^{n}(\Lambda^{\prime})=\emptyset \text{ for \ all }  1\leq n\leq MN-1.\)\\
	This is just what we would like to prove. As we shall see,
	Lemma \ref{lem-final} will  be useful in the second step.
	\medskip
	
	\textbf{Step 1.} Take $p=1, N=q, s=N_1$ and $J=\emptyset$. By Proposition \ref{prop:disp2},
there exist a positive integer $M_1$  and a $T^{M_1N}$-invariant set $\Lambda_1\subset U\cap \Lambda$
	such that
	\begin{itemize}
		\item[(A1)]  $(\Lambda_1,T^{M_1N})$ is conjugate to $(\{0,1\}^\N,\sigma)$;
		\item[(B1)]  $T^{nN}: \Lambda_1\to  T^{nN}(\Lambda_1)$ is bijective for each $1\leq n\leq M_1-1$;
		\item[(C1)] $\Lambda_1, T^N(\Lambda_1), \cdots, T^{(M_1-1)N}(\Lambda_1)$ are disjoints.
	\end{itemize}
	Furthermore, we have
	\begin{itemize}
		\item[(D1)] $\Lambda_1 \cap T^{nN+N_1}(\Lambda_1)=\emptyset  \text{ for   }   0\leq n\leq M_1-1.$
	\end{itemize}
	Indeed, on one hand,
	the fact   $U\cap  T^{N_{1}}(\Lambda)=\emptyset$ implies
	$\Lambda_1\cap T^{N_1}(\Lambda)=\emptyset$. On the other hand, the fact $\Lambda=T^{N}(\Lambda)$ implies  $T^{N_1}\Lambda=T^{nN+N_1}(\Lambda)$ for all integers $n\ge 0$.  These imply (D1) because  $\Lambda_1\subset \Lambda$.
	\medskip
	
	
	If $N=2$, we are done because both (C1) and (D1) shows that
	$\Lambda_1$ is a horseshoe with disjoint steps. If $N\ge 3$ is of the form $p^\gamma$ ($p$ being a prime), we go directly to Step 2.
	\medskip

	Using Proposition \ref{prop:disp2}, by induction on $\ell \in\{1,2,\cdots, s\}$,
	we get that there exist positive integers  $M_\ell$ ($1\le \ell\le s$)  and  a $T^{L_\ell N}$-invariant closed set $\Lambda_\ell$
	, where $L_\ell =M_1M_2\cdots M_\ell$, such that
	\begin{itemize}
		\item[(${\rm A}_\ell$)]   $(\Lambda_\ell,T^{L_\ell N})$ is conjugate to $(\{0,1\}^\N,\sigma)$;
		\item[(${\rm B}_\ell$)] $T^{nN}: \Lambda_\ell\to  T^{nN}(\Lambda_\ell)$ is bijective for $1\leq n\leq L_\ell-1$;
		\item[(${\rm C}_\ell$)] $\Lambda_\ell\cap T^{nN}(\Lambda_\ell)=\emptyset$  for   $1\leq n \leq L_s-1$;
		\item [(${\rm D}_\ell$)] $\Lambda_\ell\cap T^{nN+N_i}(\Lambda_\ell)=\emptyset$    for all $1\leq i\leq \ell$ and all $0\leq n\leq L_\ell-1$.
	\end{itemize}
	~~~~
	
	\textbf{Step 2.}
	Let  $ \{1,2, \cdots, N-1\}\setminus\{N_1, N_2, \cdots, N_s\}=\{k_1,k_2 \cdots, k_{(N-1)-s}\}$. Lemma  \ref{lem-final} says that for each $k_j$, there exist  a positive integer $n$ such  that  $n k_j\equiv N_i  (\!\!\!\!\mod N)$ for some $N_i$.
	
	Then, using Proposition \ref{prop:disp1}, by induction on
	$k_1, k_2,k_3,\cdots, k_{N-s}$
	we get an integer $M\ge 1$ and a $T^{MN}$-invariant set $\Lambda'$ such that
	\begin{itemize}
		\item[(A$^*$)]   $(\Lambda',T^{MN})$ is conjugated to $(\{0,1\}^\N,\sigma)$,
		\item[(B$^*$)] $\Lambda'\cap T^{nN}(\Lambda')=\emptyset$  for   $1\leq n \leq M-1,$
		\item [(C$^*$)] $\Lambda'\cap T^{nN+N_i}(\Lambda')=\emptyset$    for  $1\leq i\leq s$ and $0\leq n\leq M-1$,
		\item [(D$^*$)]  $\Lambda'\cap T^{nN+k_j}(\Lambda')=\emptyset$ for  $0\leq n\leq M-1$
		and $1\le j \le (N-1)-s$.
	\end{itemize}	
	Notice that the combination of (C$^*$) and (D$^*$) means
	$\Lambda'\cap T^{n}(\Lambda')=\emptyset$ for  $1\leq n\leq MN-1$, which is the desired disjointness.

\section{Proof of Theorem \ref{thm:disjiontHS}: case of two-sided Horseshoe} \label{sect:2horseshoe}

In this section, we shall prove Theorem \ref{thm:disjiontHS} for  topological dynamical systems $(X, T )$ having two-sided horseshoes $(\Lambda, T^N)$. The key is the non-existence of $n$-th root of the extended horseshoe,  stated as follows
\begin{prop}\label{prop:two-sided-no-root}
	Suppose that $(X,T)$ is  a topological dynamical system admitting a $N$-order two-sided horseshoe $(\Lambda, T^{N})$  with disjoint steps. Let $\Lambda^* =\bigcup_{k=0}^{N-1}T^k(\Lambda)$.
	Then for any integer  $n\geq 2$, there is no continuous map   $S:\Lambda^*\to \Lambda^*$ such that $S^n=T$.
\end{prop}
The proof of this proposition is given in Appendix A.
The other parts of proof are almost the same as in the case of one-sided horseshoes. We just sketch a quick proof
 by pointing out the differences from the case of one-sided horseshoe.

Notice that the maps $T^{n}  \, (n\geq 1)$ are now  bijective on $\Lambda$. This brings to something easier  because the condition (ii) concerning the bijective property of $T$  in Proposition  \ref{prop:no-root} and \ref{PO2}
is automatically satisfied.
\medskip

About the displacement of horseshoe, we have the following result.
Let   $(\Lambda,T^N)$ be a $N$-order two-sided horseshoe. Recall that,
if $\Lambda \not= T^j(\Lambda)$ for some $1\le j <N$, we say that
$\Lambda$ is displaced by $T^j$. It is true that the horseshoe $\Lambda$ is displaced by $T^j$ when $j|N$. This is actually a special case of the following Proposition \ref{prop:two-side-dis} ($p=1, q=N$ ), whose proof is
based on  Proposition \ref{prop:two-sided-no-root}.

\begin{prop}\label{prop:two-side-dis}
	Let  $(X, T )$  be  a topological dynamical system. Suppose that  $(X, T)$ has a  $pq$-order two-sided horseshoe $(\Lambda, T^{pq} )$ for some integers $p\ge 1$ and $q\ge 1$ such that
	$\Lambda\cap T^{kq}(\Lambda)=\emptyset$ for $1\leq k\leq p-1$.
	Then for any integer $1\le s < q$ such that  $s\mid q$ we have $ \Lambda\not\subset \bigcup_{k=0}^{p-1}T^{kq+s}(\Lambda)$.
\end{prop}

We omit the proof of  Proposition \ref{prop:two-side-dis}, since it is the same as  the proof of Proposition \ref{PO2}.
Note that in the present case,  $T^{n}:\Lambda \to T^{n}(\Lambda) $ is bijective for every $n\geq 1$. Hence, the condition (ii) in Proposition \ref{PO2} is automatically satisfied.
%
%
%
%
%

Taking $p=1$ and $q=N$, we obtain the following immediate corollary. In particular, $\Lambda\not \subset T(\Lambda)$.
\begin{cor}\label{PO14}
	Let $(X, T )$ be a topological dynamical system admitting a  $N$-order horseshoe $(\Lambda, T^N )$ with $N\geq 2$. Then for any positive integer $1\le s< N$ dividing $N$ we have $\Lambda\not \subset  T^s(\Lambda)$.
\end{cor}

{\it Proof of Theorem \ref{thm:disjiontHS}: case of two-sided horseshoe.}
Based on
Proposition \ref{prop:two-side-dis} (displacement of horseshoe) and   Proposition \ref{prop-SBSYS} (construction of horseshoe in a cylinder), the proof  is the same as
in the case of one-sided horseshoe. There is no need to repeat it.

\section{Appendix A. Extended horseshoes haves no roots}

The extended horseshoe has no root. A self-contained proof of this fact will be given below for the case of one-sided horseshoe. A different proof is given for the case of
two-sided horseshoe and it replies on the no existence of root for the two-sided full shift, which is known (\cite{Hedlund1969}).
It seems that the first proof can not be adapted to the case of two-sided horseshoe. 

\subsection{No root of $T: \Lambda^*\to \Lambda^*$ (one-sided case)}\ \

Here we give a proof of Proposition \ref{prop:no-root}.
The proof   is based on the investigation of the maximal entropy measure of the dynamical system  $(\Lambda^*,T)$. Given  a topological dynamical system   $(X, T)$ and an integer $s\ge 2$. A $T$-invariant measure is $T^s$-invariant. If such a measure is $T^s$-ergodic, it is $T$-ergodic. But in general, a $T$-ergodic measure is not necessarily $T^s$-ergodic. That is the case for the invariant  measure $\frac{1}{2}(\delta_{\frac{1}{3}} + \delta_{\frac{2}{3}})$ of the doubling dynamics $x \to 2 x \mod 1$. The proof of  Proposition \ref{prop:no-root} is based on the  following lemma, which shows  that the existence of  both $T$-invariant and $T^s$-ergodic measure is an obstruction
for $T^{ns}$ to be both $2$-to-$1$ and surjective.
\medskip

\begin{lem}\label{propNonexist}
	Let $(X, T)$ be a topological dynamical system
	and $s\ge 2$ an integer. Suppose that
	there exists a $T$-invariant measure which is $T^s$-ergodic. Then  $T^{ns}$ can not be both exactly $2$-to-$1$
	and  surjective for any integer $n\ge 1$.
\end{lem}

\begin{proof}
	Suppose that  $T^{ns}:X\to X$ is $2$-to-$1$ and surjective  for some integer $n\ge 1$. Then the map $T:X\to X$ must be surjective and  each point $x\in X$ has at most two pre-images by $T$, i.e. $1\leq \#T^{-1}(x) \leq 2$. This allows us to decompose $X$ into two disjoint sets
	\[A_0=\{x\in X: \#T^{-1}(x)=1\}, \qquad
	B_0=\{x\in X: \#T^{-1}(x)=2\}.\]
	
We claim that these two sets are measurable. For any $\epsilon >0$, let
	\[B_0^{\epsilon}=\{x\in X: \#T^{-1}(x)=2 \text{ and }  |y_1-y_2| \geq \epsilon
	\text{ for  distinct } y_1, y_2 \in T^{-1}(x)\}. \]
	Then $B_0^\epsilon$ is closed. Since $B_0=\cup_{n\geq 1} B_0^{1/n}$, it is an $F_\sigma$ set. Hence,  $B_0$ is measurable which implies $A_0$ is also measurable.
	

	For $0\leq k\leq ns$, let
	$B_{k}=T^{-k}(B_0)$. First notice that we have
	$$
	T^{-1}(B_k) = B_{k+1}, \quad T(B_{k+1}) = B_k \qquad (0\le k\le ns -1)
	$$
	where the second equality is because of the surjectivity of $T$.
	Since $T^{ns}$ is $2$-to-$1$, we claim that the maps $T:B_{k+1}\to B_k$ are
	injective  and  then  bijective for $1\leq k\leq ns-1$.
	
	Indeed, otherwise, for some $k$ and some $v\in B_k$ there are two distinct points  $u', u''\in B_{k+1}$  such that $T(u') = T(u'')=v$. Let $w=T^{k-1}(v)$, which belongs to $B_1=T^{-1}(B_0)$. Let $z=T(w)$, which belongs to $B_0$. By the definition of $B_0$, there exists a point $w^*\in B_1$ different from $w$ such that $T(w^*)=z$ and then there exists a point $u^*\in B_{k+1}$ different from $u', u''$ such that $T^k(u^*) = w^*$. Therefore, $z$ has at least three $T^{k+1}$-preimages $u',u'', u^*$ and then at least three $T^{ns}$-preimages. This contradicts the fact that $T^{ns}$ is $2$-to-$1$.

	The above claim implies that $B_k\subset A_0$  for $1\leq k \leq ns-1$. Since
	$A_0$ and $B_0$ form a partition of $X$, we have   $B_{k}\cap B_0=\emptyset$ for $1\leq k \leq ns-1$. Consequently, all  $B_j$'s for
	$0\leq j \leq ns-1$ are disjoint.
	We claim that all  $B_j$'s for
	$0\leq j \leq ns-1$ form a partition of $X$. To that end, it suffices to prove
	\begin{align}\label{eqA0}
		A_0=\bigsqcup_{k=1}^{ns-1}B_{k}.
	\end{align}
	Suppose that (\ref{eqA0}) is not true, which means there exists a point $x\in A_0\setminus\bigsqcup_{k=1}^{ns-1}B_{k}$. Hence $T^{k}(x)\in A_0$ for $0\leq k\leq ns-1 $, by the definition of $B_k$ and the fact that $\{A_0, B_0\}$ is a partition of $X$. Therefore, the point $T^{ns-1}(x)\in A_0$ has a unique $T^{ns}$-preimage $T^{-1}(x)$, which contradicts to the assumption that $T^{ns}$ is $2$-to-$1$.

	Since $X= A_0\sqcup B_0$, we get the decomposition
	\[X=\bigsqcup_{k=0}^{ns-1}B_{k}
	= \bigsqcup_{j=0}^{s-1} T^{-j} (X') \ \ \
	{\rm with }\ \ X^{\prime}=\bigsqcup_{k=0}^{n-1}B_{ks}.
	\]
	By the hypothesis, there exists a  $T$-invariant measure $\mu$ which is  $T^s$-ergodic.
	The $T$-invariance of $\mu$ implies that $\mu(X^{\prime})=1/s$.

	Let $A_1=T^{-1}(A_0)$. Recall that $B_1 =T^{-1}(B_0)$.
	Since $\{A_0, B_0\}$ is a partition of $X$, so is
	$\{A_1, B_1\}$. As $ B_1\subset A_0$ which is proved above,  we have $B_0\subset A_1$ and actually  $B_0=A_1\setminus A_0$.    By \eqref{eqA0} and the definition of $A_1$, we get $A_1= \bigsqcup_{k=2}^{ns}B_{k}$. Then from
	$A_0\sqcup B_0 = A_1\sqcup B_1$ we get  $B_0=B_{ns}$, in other words,
	$$
	T^{-ns}(B_0)=B_0.$$
	Consequently,  $T^{-s}(X^{\prime})=X^{\prime}$. Then, by the $T^s$-ergodicity of $\mu$, we have $\mu(X^{\prime})=0 \text{ or } 1$, which  contradicts to $\mu(X') =\frac{1}{s}$ with $s\ge 2$.
\end{proof}

Now we shall prove Proposition \ref{prop:no-root} by contradiction. Basic properties of entropy function will be used. We refer to Walters' book \cite{Walters1982} (cf. Theorem 4.13, Theorem 7.5, Theorem 7.10, Theorem 8.1).

\begin{proof}[Proof of Proposition \ref{prop:no-root}]
	Assume that   $S:\Lambda^*\to \Lambda^*$  is   a  continuous map such that $S^n=T$ for some integer $n\geq 2$.

	For $1\le j\le N-1$, we  consider the dynamical system  $(T^{j}(\Lambda), T^{N})$.
	By the hypothesis,  the map $T^j:\Lambda\to T^{j}(\Lambda)$ is a homeomorphism. So,  the system $(T^{j}(\Lambda), T^{N})$  is conjugate to $(\Lambda, T^{N})$ with the conjugation $T^{j}$. Hence, the topological entropy of the system $(T^{j}(\Lambda), T^{N})$ is equal to
	$ \log 2$. It follows that
	$
	h_{\rm top}(\Lambda^*,T^{N})=\log 2$. Therefore,
	$$
	h_{\rm top}(\Lambda^*,T)=\frac{1}{N}\log 2.$$
	It follows that the factor
	dynamical  system $(\Lambda^*,S)$ admits its topological entropy
	\[h_{\rm top}(\Lambda^*,S)=\frac{1}{nN}\log 2.\]
	Let $\mu$ be  a maximal entropy  measure of the system $(\Lambda^*,S)$.   It  is also a maximal entropy measure of $(\Lambda^*, T)$. By the same reason,  $\mu$ is also a maximal entropy measure of dynamical system $(\Lambda^*,T^{N})$.
	
	Let $\mu_j = \mu|_{T^{j}(\Lambda)}$ be the restriction
	for $0\le j\le N-1$. By the $T$-invariance of $\mu$ and the disjointness of $T^{j}(\Lambda)$'s ($0\le j \le N-1$), it is easy to get
	$$
	\mu_{j+1} = \mu_j \circ T^{-1}   \quad (0\le j \le N-1)
	$$
	with the convention $\mu_N =\mu_0$. It follows that $\mu_j$ are all $T^{N}$-invariant.
	Let $\nu_j= N\cdot\mu_j$, which is a probability measure concentrated on $T^{j}(\Lambda)$. Since the systems $(T^{j}(\Lambda), \nu_j, T^{N})$ are all conjugate, the measure-theoretic entropies $h_{\nu_j}(T^{j}(\Lambda), T^{N})$
	are equal and the common value is $h_{\nu_j}(\Lambda^*, T^{N})$. From this, the relation $\mu= \frac{1}{N} \sum_{j=0}^{N-1}\nu_j$ and the
	affinity of entropy, we get
	\[
	\log 2=h_\mu(\Lambda^*, T^{N})=
	\frac{1}{N}\sum_{j=0}^{N-1} h_{\nu_j}(\Lambda^*, T^{N})
	= h_{\nu_0}(\Lambda^*, T^{N}) = h_{\nu_0}(\Lambda, T^{N}).
	\]
	So, the measure $\nu_0$ is  the unique maximal entropy measure of the horseshoe $(\Lambda, T^{N})$,
	i.e. the  Bernoulli $(\frac{1}{2},\frac{1}{2})$-measure on $\{0,1\}^{\N}$. Since $\nu_0$ is $T^{N}$-ergodic,   so are $\nu_j$'s. Now we claim that $\mu$ is $T$-ergodic.
	In fact, assume that $A\subset \Lambda^*$ is a $T$-invariant set. Then $A\cap T^{j}(\Lambda)$ is $T^{N}$-invariant for every $0\le j\le N-1$. By the $T^{N}$-ergodicity of $\nu_j$,  $\nu_j(A\cap T^{j}(\Lambda)) =0 \ {\rm or} \ 1$ for each $j$. But $\nu_j(A\cap T^{j}(\Lambda))$ are equal for different $j$'s. Then  we get $\mu(A) =0 \ {\rm or}\ 1$, because
	$$
	\mu(A) = \frac{1}{N}\sum_{j=0}^{N-1} \nu_j(A\cap T^{j}(\Lambda)).
	$$

The existence of $S$-invariant measure $\mu$ which is  $T$-ergodic measure,  contradicts the fact
	that $T^N: \Lambda^*\to \Lambda^*$ is $2$-to-$1$, by Lemma \ref{propNonexist}.
\end{proof}

\subsection{No root of $T: \Lambda^*\to \Lambda^*$ (two-sided case)}
\label{sect:noroot}\ \

Here we prove Proposition \ref{prop:two-sided-no-root}.
The proof is based on the fact that the shift map $\sigma: \{0,1\}^\mathbb{Z}\to \{0,1\}^\mathbb{Z}$ has no root, which is well known (cf. \cite{Hedlund1969}, Corollary 18.2, p.371).

Given a  two-sided horseshoe $(\Lambda,T^{N})$ with disjoint steps, it is convenient to identify the subsystem  $(\Lambda^*,T)$  with the following system $(\{0,1\}^{\Z}\times \Z/N\Z, \sigma_N)$, a tower of height $N$, where $\sigma_N$ is defined by
	\begin{align}\label{def:sigman}
\sigma_N(\omega,k)=\left\{
	\begin{array}{ll}
	(\sigma(\omega),0), & \mbox{if $k=N-1$;} \\
	(\omega,k+1), & \mbox{otherwise.}
	\end{array}
	\right.
	\end{align}
	The tower $\{0,1\}^{\Z}\times \Z/N\Z$ has $N$ floors $F_i = \{0,1\}^{\Z}\times \{i\}$ for $0\le i <N$. We extend
	$F_i$ for all integers $i\ge 0$ by defining $F_i = F_{i \mod N}$. Especially $F_N = F_0$.
	In the following we denote $\sigma_N$ by $T$.

Suppose that $T$ has a root $S$, i.e. $S^p=T$ for some $p\ge 2$. It is clear that $S$ is bijective and commutes with $\sigma_N$.
We claim that $S$ permute floors, that is to say, for any $i$ there exists a $j$ such that $S (F_i) =F_j$. We shall prove this claim  by the commutativity of $S$ with $T$.  Assume this claim for the moment. Then each floor $F_i$  ($0\le i <N$) is mapped back to $F_i$ by $S^N$,
that is to say,
\begin{equation}
 \forall \omega, \quad S^N(\omega, i) = (R_i\omega, i),
\end{equation}
where $R_i: \{0, 1\}^\mathbb{Z}\to \{0, 1\}^\mathbb{Z}$ is some map, which depends on $i$. It is easy to see that $R_i$ is continuous. Thus, on one hand, we have
$$
       \forall \omega, \quad S^{Np}(\omega, i) =(R_i^p\omega, i);
$$
on the other hand, we have
$$
       \forall \omega, \quad S^{Np}(\omega, i) = T^N( \omega, i)= ( \sigma \omega, i).
$$
It follows that $\sigma=R_i^p$, which is impossible.

Now let us prove the claim. Let $\mu$ the ergodic probability measure of the maximal entropy of $T$ (its restriction on each floor is  the symmetric Bernoulli measure).
We have $\mu(F_i) =\frac{1}{N}$.  Let  $C_i = S(F_0)\cap F_i$, the portion of $SF_0$ contained in $F_i$, for $0\le i <N$.  Suppose $0<\mu(C_i)<\mu(F_i)$
for some $i$. There would be a contradiction. Indeed, we have the invariance
$$
   T \left(\bigcup_{j=0}^{N-1} T^j (C_i)\right) = \bigcup_{j=0}^{N-1} T^j (C_i).
$$
This is because
$$
     T^j(C_i)=  T^j(S(F_0)\cap F_i)= ST^j (F_0) \cap T^jF_{i} = S(F_j) \cap F_{i+j}
$$
and
$$
   \bigcup_{j=1}^{N}  S(F_j) \cap F_{i+j} = \bigcup_{j=0}^{N-1}  S(F_j) \cap F_{i+j}.
$$
The invariant set $\bigcup_{j=0}^{N-1} T^j (C_i)$ has its measure between $0$ and $1$ because of $0<\mu(C_i)<\frac{1}{N}$. This contradicts the ergodicity of $\mu$.
We have thus proved that for every $i$, the measure $\mu(S (F_0) \cap F_i)$ is equal to $0$ or $\mu(F_i)$.  There exists one $i$ such that $\mu(S (F_0) \cap F_i)= \mu(F_i)$, otherwise
$\mu(S(F_0))=0$ so that
$$\mu(S(F_j))= \mu(ST^j (F_0))= \mu(T^j S (F_0)) = \mu(S(F_0))=0,$$
which implies that $\mu$ is the null measure. There is at most one $i$ such that $\mu(S (F_0) \cap F_i)= \mu(F_i)$, otherwise $\mu(S(F_0))\ge \frac{2}{N}$, which implies that $\mu$
has a total measure equal to at least $2$.  So, $S(F_0)$ and $F_i$ are equal almost everywhere. If we take into account the continuity of $S$, we get that the two compact sets
$S(F_0)$ and $F_i$ must be
equal. In place of $F_0$, we can consider any $F_j$. The same argument shows that $S(F_j)$ must be equal to some $F_k$.  In this way, $S$
defines a permutation on floors. Otherwise, under $S$ we have a cycle
     $$
       F_{i_0} \to F_{i_1}\to \cdots \to F_{i_{\ell-1}}\to F_{i_0}
     $$
     with $\ell <N$. Then the union $U$ of these $\ell$ floors is $S$-invariant and it  is also $S^p$-invariant, i.e. $T$-invariant, an obvious contradiction.

\section{Appendix B.  Any cylinder contains a horseshoe}

\subsection{Horseshoe with disjoint steps in any cylinder: one-sided case}\ \

Consider the one-sided full shift system $(\{0,1\}^{\N}, \sigma)$. For any non-empty open set  $U\subset \{0,1\}^{\N}$, we shall show that there exists a horseshoe $(\Lambda,\sigma^N)$ with disjoint steps for some sufficiently large integer $N\geq 1$ such that
$ \Lambda\subset U$ and   the map   $\sigma^{N-1}: \Lambda\to \sigma^{N-1}(\Lambda)$ is bijective. In other word, this horseshoe $(\Lambda,\sigma^N)$ satisfies the conditions (i) and (ii) required by Proposition \ref{prop:no-root}.

For a word $a_0\cdots a_{n-1}$ of length $n$, denote  by   $[a_0\cdots a_n-1]$ the cylinder of rank $n$
\[[a_0\cdots a_{n-1}]=\{y \in  \{0,1\}^{\N}: y_i=a_i \text{ for } 0\leq i \leq n-1\}.\]
Denote the $n$-prefix of $x=(x_j)_{j\ge 0}$ by $x|_n$, i.e. $x|_n=x_0x_1\cdots x_{n-1}$.
Hence $[x|_{n}]$ denote a cylinder of rank $n$.
Let $uv$ denote the concatenation $u_0u_1\cdots u_{n-1}v_0\cdots v_{m-1}$ of two words $u=u_0u_1\cdots u_{n-1}$ and  $v= v_0\cdots v_{m-1}$. So, $u^{r}$ denotes $u\cdots u$ ($r$ times). In particular, $1^r$ means $1\cdots 1$ ($r$ times).  The following lemma is a preparation for proving the above announced existence of horseshoe in a given cylinder.

\begin{lem} \label{lem-EIF}
	For any cylinder $C$ in
	$\{0,1\}^\N$ of rank $M\ge 1$, there exists a sub-cylinder $C'\subset C$ of rank $N$ with $N\geq M$   such that
	\begin{itemize}
		\item[(i)]$\sigma^N (C') = \{0,1\}^\N; $
		\item[(ii)]$C'\cap \sigma^{n}(C')=\emptyset \ \text{for } 1\leq n \leq N-1.$
	\end{itemize}
	
\end{lem}

\begin{proof}
	Assume $C=[a_0a_1\cdots a_{M-1}]$. We assume $a_0=0$, without loss of generality. R Note that $\sigma^{M}(C)=\{0,1\}^{\N}$. Let  $n_0\ge 1 $ be the minimal positive integer such that  $C \subset \sigma^{n_0}(C)$, so that
	\begin{align}\label{eq:CCapC}
		C\cap \sigma^{n}(C)=\emptyset \quad\  \text{for\ all }\ 1\leq n\leq n_0-1.
	\end{align}
	We restate these facts as follow
	\begin{align}\label{eq:disj1}
		[0a_1\cdots a_{M-1}]\cap  [a_1\cdots a_{M-1}]=\emptyset ,& \cdots,  [0a_1\cdots a_{M-1}]\cap [a_{n_0-1}a_{n_0}\cdots a_{M-1}]=\emptyset;
	\end{align}
	\begin{align}\label{eq:includsion}
		[0a_1\cdots a_{M-1}]\subset [a_{n_0}\cdots a_{M-1}].
	\end{align}
	We have $n_0\le M$. If $n_0=M$, we are done and we can take $C'=C$. In the following, we assume $n_0<M$.
	
	The inclusion \eqref{eq:includsion} means
	\begin{align}\label{eq:perodic}
		a_{n_0+j}=a_j  \quad  \text{for } 0\leq j \leq M-n_0-1.
	\end{align}
	Let $x=\overline{a_0a_1\cdots a_{n_0-1}}$, which is $n_0$-periodic. By \eqref{eq:perodic}, $a_0a_1\cdots a_{M-1}$ is a prefix of $x$, so $x$ is in $C$.
	By (\ref{eq:disj1}), $n_0$ is the minimal period of $x$. Define
	the sub-cylinder $C^{\prime}=[0a_1a_2\cdots a_{M-1}1^{n_0}]$, or more precisely
	$$
	C^{\prime} =[(0a_1\cdots a_{n_0-1})^{q}0a_1\cdots a_{j-1}1^{n_0}]
	$$
	where $q\ge 0$ and $0\leq j \leq n_0-1$ are determined by $M=qn_0+j$.  Now we shall check that the sub-cylinder $C^{\prime}$(of $C$) of rank $N=M+n_0$  has property (ii). We distinguish three cases.
	
	{\it Case I. $1\leq n\leq  n_{0}-1$}. Since $C^{\prime}\subset C$,  we proved $C^{\prime}\cap \sigma^{n}(C^{\prime})=\emptyset$, by \eqref{eq:CCapC}.
	
	{\it Case II. $n_0\leq n \leq M-1$}. Suppose  $C^{\prime}\cap \sigma^{n}(C^{\prime})\neq \emptyset$. Then
	$C\cap \sigma^{n}(C^{\prime})\neq \emptyset$. Since $|\sigma^{n}(C^{\prime})|\leq M=|C|$, we get
	$C\subset \sigma^{n}(C^{\prime})$.  Hence, the word $a_n\cdots a_{M-1}1^{n_0}$ defining the cylinder $\sigma^{n}(C^{\prime})$ has $1^{n_0}$ as suffix, which is a word contained in $a_0a_1\cdots a_{M-1}$ defining the cylinder $C$. But, on the other hand,  any word of length $n_0$ contained in $a_0a_1\cdots a_{M-1}$ contains $0$, a contradiction.

	{\it Case III. $M\leq n \leq N-1$}. This case is evident because $\sigma^{n}(C^{\prime})\subset [1]$, but $C'\subset [0]$.
\end{proof}

Let us look at the cylinders $C$ of rank $M=4$ and the cylinders $C'$ constructed in Lemma \ref{lem-EIF}:
\begin{eqnarray*}
	C=[0000], & C' = [00001^1]\\
	C=[0001], & C' = C\ \ \ \ \ \ \ \ \\
	C=[0010], &  C' = [00101^3]\\
	C=[0011], & C' = C\ \ \ \ \ \ \ \ \\
	C=[0100], & C' = [01001^3]\\	
	C=[0101], & C' = [01011^2]\\
	C=[0110], & C' = [01101^3]\\
	C=[0111], &   C' = C\ \ \ \ \ \ \ \
\end{eqnarray*}
where the exponent represents $n_0$.

We are now ready to prove the existence of horseshoe contained in a given cylinder.

\begin{prop} \label{lem-SBSYS}
	Let $C\subset\{0,1\}^{\N}$ be an arbitrary cylinder.  There exists  a $\sigma^{N}$-invariant  closed subset $\Lambda\subset C$ for some integer $N\geq |C|$, such that
	\begin{itemize}
		\item[(i)] The system $(\Lambda, \sigma^N)$ is topologically conjugate to the full shift $(\{0,1\}^{\N}, \sigma)$;
		\item[(ii)] The maps $\sigma^{n}: \Lambda \to  \sigma^{n}(\Lambda) \ (1\leq n\le N-1)$ are bijections;
		\item[(iii)] The sets $\Lambda,  \sigma (\Lambda), \cdots, \sigma^{N-1}(\Lambda)$ are disjoint.
	\end{itemize}
\end{prop}
\begin{proof}
	Assume   $C=[y_0y_1\cdots y_{m-1}]$ be a  cylinder of rank $m$. We can assume $y_0=0$ without loss of generality.
	By Lemma \ref{lem-EIF}, there exists  an integer $n_* \geq  m$ and a sub-cylinder $C^{\prime}\subset C$ of rank $n_*$ such that
	\begin{equation}\label{eq:disj}
	C^{\prime}\cap \sigma^{n}(C^{\prime})=\emptyset\ \  \text{ for } 1\leq n \leq n_*-1.
	\end{equation}
	It is obvious that $\sigma^{n_*}(C^{\prime})=\{0,1\}^\mathbb{N}$
	and  $\sigma^{n_*}: C' \to \{0,1\}^\mathbb{N}$ is  bijective.
	Assume that $C' =[y_0y_1\cdots y_{m-1} y_m \cdots y_{n_*-1}]$.
	Notice that $y_{n_*-1}=1$ because $C'\cap  \sigma^{n_*-1}(C')=\emptyset$.
	Let $y= \overline{y_0y_1\cdots y_{n_*-1}}$. This periodic point $y$ is the unique periodic point contained in $C^{\prime}$ of exact period $n_*$, by (\ref{eq:disj}).
	Let $n_1=\min\{0\le i\le n_*-1: y_i=1\}$.  Since $y_0=0$ and $y_{n_*-1}=1$,
	we have
	$1\leq n_1\leq n_*-1$.
	Define   two sub-cylinders of $C'$ of rank $n_*+n_1+2$:
	\begin{eqnarray*}
		C_1 &= &[y_0y_1\cdots y_{n_*-1}10  0^{n_1}]
		=[0^{n_1}y_{n_1}\cdots y_{n_*-1}10  0^{n_1}]=[0^{n_1}a0^{n_1}],
		\\
		C_2 &= &[y_0y_1\cdots y_{n_*-1}11 0^{n_1}] =
		[0^{n_1}y_{n_1}\cdots y_{n_*-1}11 0^{n_1}] =[0^{n_1}b0^{n_1}]
	\end{eqnarray*}
	where $a = y_{n_1}\cdots y_{n_*-1}10$
	and $b = y_{n_1}\cdots y_{n_*-1}11$.

	Also observe that \\
	\indent (iv)
	{\it for $1\leq n\leq n_*-1$, $\sigma^{n}(C_1\cup C_2)\cap C^{\prime}=\emptyset $  and  $\sigma^{n}$ is injective on $C_1\cup C_2$}.\\
	This follows from the relation $C_{1}\cup C_2 \subset C^{\prime}$, the disjointness (\ref{eq:disj})
	and
	the injectivity of  $\sigma^{n_*}: C^{\prime} \to \{0,1\}^{\N}$.    By the definition of $n_1$,  we have
	$C_{1}\cup C_{2}\subset [0^{n_1}1]$.
	Note that  \[\sigma^{n_*}(C_1)=[100^{n_1}], \quad \sigma^{n_*+1}(C_1)=[0^{n_1+1}]; \quad
	\sigma^{n_*}(C_2)=[110^{n_1}], \quad \sigma^{n_*+1}(C_2)=[10^{n_1}]\]
	which are all disjoint from $[0^{n_1}1]$.
	Hence,   \\
	\indent (v) {\it for  $ n=n_*$ or $n_*+1$,
		$\sigma^{n}(C_1\cup C_2)\cap (C_{1}\cup C_{2})=\emptyset $
		and  $\sigma^n$ is injective on $C_1\cup C_2$}.

	Now take $N=n_*+2$ and  define
	\[
	\Lambda=\{x\in\{0,1\}^{\N}: \forall k\ge 0, \sigma^{kN}(x) \in  C_{1}\cup C_{2}\}.
	\]
	Observe that both the words defining $C_1$
	and $C_2$ have $0^{n_1}$ as their prefix and as well as their suffix. Therefore $\Lambda$ can be identified with the symbolic space $\{0^{n_1}a, 0^{n_1} b\}^\mathbb{N}$ and $(\Lambda, \sigma^N)$ is conjugate to the shift map on $\{0^{n_1}a, 0^{n_1} b\}^\mathbb{N}$. This is the property (i) of $(\Lambda, \sigma^N)$.
	The required properties (ii) and (iii) follow from (iv) and (v).
	%
	%
\end{proof}

\subsection{Horseshoe with disjoint steps in any cylinder: two-sided case}

Consider the full shift system $(\{0,1\}^{\Z}, \sigma)$. For any non-empty cylinder  $C\subset \{0,1\}^{\Z}$, we shall show that there exists a horseshoe $(\Lambda,\sigma^N)$ with disjoint steps for some sufficiently large integer $N\geq 1$, where $\Lambda \subset C$. We start with the following preparative lemma, the counterpart of Lemma \ref{lem-EIF}.

\begin{lem} \label{lem-EIF4}
	For any cylinder $C$ in
	$\{0,1\}^\Z$ of rank $M\ge 1$, there exists a sub-cylinder $C'\subset C$ of rank $N$ with $N\geq M$   such that
	\begin{itemize}
		\item[(i)]$\sigma^N (C') \cap  C'\neq \emptyset; $
		\item[(ii)]$C'\cap \sigma^{n}(C')=\emptyset \ \text{for } 1\leq n \leq N-1.$
	\end{itemize}
	
\end{lem}

\begin{proof}

	The proof is almost the same as that of Lemma \ref{lem-EIF}.
	We assume $C=[a_0a_1\cdots a_{M-1}]$,  and $a_0=0$ without loss of generality. Note that $C\cap \sigma^{M}(C)\neq \emptyset $. Let  $n_0\ge 1 $ be the minimal positive integer such that  $C\cap \sigma^{n_0}(C)\neq \emptyset $, so that
	\begin{align}\label{eq:CCapC4}
		C\cap \sigma^{n}(C)=\emptyset \quad\  \text{for\ all }\ 1\leq n\leq n_0-1.
	\end{align}
	We restate these facts as follow
	\begin{align}\label{eq:disj14}
		[0a_1\cdots a_{M-1}]_0\cap  [a_1\cdots a_{M-1}]_0=\emptyset ,& \cdots,  [0a_1\cdots a_{M-1}]_0\cap [a_{n_0-1}a_{n_0}\cdots a_{M-1}]_0=\emptyset;
	\end{align}
	\begin{align}\label{eq:includsion4}
		[0a_1\cdots a_{M-1}]_0\subset [a_{n_0}\cdots a_{M-1}]_0.
	\end{align}
The rest of the proof is the same if we replace cylinders of the for $[*]$ by $[*]_0$ where $*$ represents a word.
\end{proof}

We are now ready to prove the existence of horseshoe contained in a given cylinder.

\begin{prop} \label{prop-SBSYS}
	Let $C\subset\{0,1\}^{\N}$ be an arbitrary cylinder.  There exists  a $\sigma^{N}$-invariant  closed subset $\Lambda\subset C$ for some integer $N\geq|C|$, such that
	\begin{itemize}
		\item[(i)] The system $(\Lambda, \sigma^N)$ is topologically conjugate to the full shift $(\{0,1\}^{\N}, \sigma)$;
		
		\item[(ii)] The sets $\Lambda,  \sigma (\Lambda), \cdots, \sigma^{N-1}(\Lambda)$ are disjoint.
	\end{itemize}
\end{prop}
\begin{proof}
	It is exactly the same proof as that of Proposition \ref{lem-SBSYS},
	if we replace cylinders of the form $[*]$ by $[*]_0$.
\end{proof}

\end{document}